\documentstyle[11pt]{article}
\begin{document}
\newtheorem{prop}{Proposition}[section]
\newtheorem{Def}{Definition}[section]
\newtheorem{theorem}{Theorem}[section]
\newtheorem{lemma}{Lemma}[section]
\newtheorem{Cor}{Corollary}[section]

\title{\bf Global solutions of the Klein - Gordon - \\
Schr\"odinger system with rough data}
\author{{\bf Hartmut Pecher}\\Fachbereich Mathematik\\ Bergische Universit\"at Wuppertal\\ Gau{\ss}str. 20 \\ D-42097 Wuppertal\\ Germany\\ e-mail Hartmut.Pecher@math.uni-wuppertal.de}
\date{}
\maketitle

\begin{abstract}
The Klein - Gordon -- Schr\"odinger system with Yukawa coupling is shown to have a unique global solution for rough data, which not necessarily have finite energy. The proof uses a generalized bilinear estimate of Strichartz type and Bourgain's idea to split the data into low and high frequency parts.
\end{abstract}

\normalsize
\setcounter{section}{-1}
\section{Introduction}
In this paper we consider the Cauchy problem for the (3+1)-dimensional \\ Klein - Gordon - Schr\"odinger system with Yukawa coupling:
\begin{eqnarray}
\label{-1.1}
i\psi_t + \Delta \psi & = & -\phi\psi \\
\label{-1.2}
\phi_{tt} + (-\Delta + 1)\phi & = & |\psi|^2 \\
\label{-1.3}
\psi(0) \; = \; \psi_0 \quad , \quad \phi(0) & = & \phi_0 \quad , \quad \phi_t(0)\; = \;\phi_1
\end{eqnarray}
Here $ \psi: {\bf R}^3 \times {\bf R} \to {\bf C} $ is the nucleon field and $ \phi: {\bf R}^3 \times {\bf R} \to {\bf R}$ the meson field.
This system satisfies the conservation laws (\ref{2.1}) and (\ref{2.2}) below. It is well-known that these conservation laws imply the existence of a unique global solution for data $ \psi_0 \in H^{1,2}({\bf R}^3) , \phi_0 \in H^{1,2}({\bf R}^3) , \phi_1 \in L^2({\bf R}^3) $ with $ \psi \in C^0({\bf R},H^{1,2}({\bf R}^3)) , \phi \in C^0({\bf R},H^{1,2}({\bf R}^3)) \cap C^1({\bf R},L^2({\bf R}^3)) $ (cf. \cite{B}, Thme. 3). The same is true for data $\psi_0 \in H^{m,2}({\bf R}^3) , \phi_0 \in H^{m,2}({\bf R}^3) , \phi_1 \in H^{m-1,2}({\bf R}^3) , m \in {\bf N} , m \ge 2 $ (cf. \cite{BC}).

Our main results are local existence and uniqueness for data $ (\psi_0,\phi_0,\phi_1) \in H^{s,2}({\bf R}^3) \times H^{m,2}({\bf R}^3) \times H^{m-1,2}({\bf R}^3) $ with $s=m=0$ and more regular data (Theorem \ref{Theorem 1.3}) and global existence and uniqueness for data without finite energy, namely under the assumptions
$ 1 \ge s,m > 7/10 $ and $ s+m > 3/2 $ (Theorem \ref{Theorem 5.1}).
 
We use Bourgain's idea (\cite{B1},\cite{B2}) to split the data into low and high frequency parts in order to prove global well-posedness for rough data for which the conservation laws are not directly applicable. This technique has been successfully applied to other problems (\cite{B5},\cite{CST},\cite{FLP},\cite{KT},\cite{KPV1},\cite{P}). We should also mention\\[0.3cm] 
$\overline{\scriptstyle { 2000 \,  Mathem}}{\scriptstyle atics \, Subject\, Classification \,35Q99}$\\
that for problems with a scaling invariance (which is not true in our situation) improved and sometimes optimal global-wellposedness results were 
given by Colliander-Keel-Staffilani-Takaoka-Tao using the so-called ''I-method`` (\cite{CKSTT1},\cite{CKSTT2}, \cite{CKSTT3},\cite{CKSTT4},\cite{CKSTT5}).

We also rely on \cite{GTV} for the local theory and the framework for the technique. The key point is a new smooth
ing estimate for the nonlinearity $\phi \psi$ in equation (\ref{-1.1}) given in Lemma \ref{Lemma 0.5}
which is a generalized bilinear Strichartz type estimate. A similar estimate has been given before by Bourgain for the pure Schr\"odinger problem (\cite{B1},\cite{B2}). 

The paper is organized as follows. In Section 1 the estimates for the nonlinearities are given along the lines of Ginibre, Tsutsumi and Velo \cite{GTV} in the $X^{s,b}$-spaces, which were introduced by Bourgain (\cite{B3},\cite{B4}), Klainerman and Machedon (\cite{KM1},\cite{KM2}) and Kenig, Ponce and Vega \cite{KPV}.

For an equation of the form 
\begin{equation} 
\label{-1.4}
iu_t - \varphi(-i\nabla_x)u = 0 
\end{equation}
where $\varphi$ is a measurable function, let $X_{\varphi}^{s,b}$ be the completion of
$S({\bf R}^4)$ with respect to
\begin{eqnarray*}
\| f \|_{X^{s,b}_{\varphi}} : & = & \| <\xi>^s <\tau>^b {\cal F}(e^{it\varphi(-i\nabla_x)}f(x,t)) \|_{L^2_{\xi\tau}} \\
& = & \| <\xi>^s <\tau+\varphi(\xi)>^b \widehat{f}(\xi,\tau) \|_{L^2_{\xi,\tau}}
\end{eqnarray*}
For $ \varphi(\xi) = \pm < \xi > $ we use the notation $ X^{s,b}_{\pm} $ and for $ \varphi(\xi) = | \xi |^2 $ simply $ X^{s,b} $. For a given time interval $I$ we define 
$$ \| f \|_{X^{s,b}(I)} = \inf_{\widetilde{f}_{|I}=f} \| \widetilde{f} \|_{X^{s,b}} \quad \mbox{and similarly} \quad  \| f \|_{X^{s,b}_{\pm}(I)} $$ 
In Section 2 we transform the system into a first-order system for $(\psi,\phi_+,\phi_-)$ and give the local well-posedness result as a variation of a result of \cite{GTV} for the Zakharov system. Especially we have local solutions in all those cases where we want to show the existence of a unique global solution.

In Section 3 we split the data into sums $ \psi_0 = \psi_{01} + \psi_{02} \, , \, \phi_0 = \phi_{01} + \phi_{02} \, , \, \phi_1 = \phi_{11} + \phi_{12} $ , where the low frequency parts $(\psi_{01},\phi_{01},\phi_{11})$ are regular with large norms and the high frequency parts $(\psi_{02},\phi_{02},\phi_{12})$ are just in $H^{s,2} \times H^{m,2} \times H^{m-1,2} $ with small $L^2$-norms.

In Section 4 the solution $(\widetilde{\psi},\widetilde{\phi}_+,\widetilde{\phi}_-)$ of the corresponding first order problem with data $(\psi_{01},\phi_{0+1},\phi_{0-1})$ is further investigated on a suitable time interval $I$ using the conservation laws and Strichartz type estimates.

In Section 5 we consider the system fulfilled by $ (\widehat{\psi},\widehat{\phi}_{\pm}) = (\psi - \widetilde{\psi}, \phi_{\pm} - \widetilde{\phi}_{\pm})$ with data $(\psi_{02},\phi_{0\pm 2})$ and construct a solution on the same time interval $I$, thus giving a solution of our original problem on $I$. The inhomogeneous parts $w$ of $\widehat{\psi}$ and $z_{\pm}$ of $\widehat{\phi}_{\pm}$ are shown to belong to $H^{1,2}({\bf R}^3)$, thus to be more regular than the homogeneous parts $e^{it\Delta}\psi_{02}$ and $e^{\mp it(-\Delta +1)^{1/2}}\psi_{0\pm 2}$ which belong to $H^{s,2}({\bf R}^3)$ and $H^{m,2}({\bf R}^3)$, respectively.

In Section 6 we show that the process can be iterated to get a solution on any time interval $[0,T]$. One constructs solutions step by step on time intervals of equal length $|I|$, taking as new initial data at time $|I|$ the triple $ (\widetilde{\psi}(|I|)+w(|I|),\widetilde{\phi}_{\pm}(|I|)+z_{\pm}(|I|))$ and repeating the argument on $[|I|,2|I|]$. In each step the involved norms have to be controlled independently of the iteration step in order to be able to choose intervals of equal length.

We use the following standard facts about the spaces $X^{s,b}_{\varphi}$:\\
if $u$ is a solution of (\ref{-1.4}) with $u(0)=f$ and $\psi$ a cutoff function in $C^{\infty}_0({\bf R})$ with $ supp \, \psi \subset (-2,2) $ , $ \psi \equiv 1$ on $[-1,1] $ , $ \psi(t) = \psi(-t) $ , $ \psi(t) \ge 0 $ , $ \psi_{\delta}(t) := \psi(t/\delta) $ ,  if $ 0 < \delta \le 1 $ , we have, for $ b \ge 0 $ :
\begin{equation}
 \label{-1.4a}
\| \psi_1 u \|_{X^{s,b}_{\varphi}} \le c \| f \|_{H^{s,2}_x} 
\end{equation}
If $v$ is the solution of 
$$ iv_t + \varphi(-i\nabla_x)v = F \, , \, v(0) = 0 $$
we have, for $b'+1 \ge b \ge 0 \ge b' > -1/2 $:
\begin{equation}
 \label{-1.5}
\| \psi_{\delta}v\|_{X^{s,b}_{\varphi}} \le c \delta^{1+b'-b} \|F\|_{X^{s,b'}_{\varphi}} 
\end{equation}
(for a proof see \cite{GTV}, Lemma 2.1).\\
We have $ X^{s,b}_{\varphi}(I) \subset C^0(I,H^{s,2}({\bf R}^3)) $ , if $ b > 1/2 $ , $ I \subset {\bf R} $. \\
The following interpolation property holds:
$$ X^{s,b}_{\varphi}(I) = (X^{s_0,b_0}_{\varphi}(I),X^{s_1,b_1}_{\varphi}(I))_{[\Theta]} $$
where $ 0 \le \Theta \le 1 \, , \, b=(1-\Theta) b_0 + \Theta b_1 \, , \, s = (1-\Theta)s_0 + \Theta s_1 $ (see \cite{P}, Chapter 0).\\
We also use the following consequences of the Strichartz type estimates for the Schr\"odinger equation in ${\bf R}^n$:\\
For $ 0 \le \frac{2}{q} = n(\frac{1}{2}-\frac{1}{r}) < 1 $ the following estimate holds
\begin{equation}
 \label{-1.6}
\|e^{it\Delta}\psi_0\|_{L^q_t(I,L^r_x({\bf R}^n))} \le c \| \psi_0\|_{L^2_x({\bf R}^n)} 
\end{equation}
and implies
\begin{equation}
 \label{-1.7}
\| f \|_{L^q_t(I,L^r_x({\bf R}^n))} \le c \| f \|_{X^{0,\frac{1}{2}+}(I)} 
\end{equation}
The dual version also holds:
\begin{equation}
 \label{-1.8}
\| f \|_{X^{0,-\frac{1}{2}-}(I)} \le c \| f \|_{L^{q'}_t(I,L^{r'}_x({\bf R}^n))}
\end{equation}
where $q',r'$ are dual to $q,r$.\\
Similarly, for the Klein-Gordon equation, let $ 0 \le \frac{2}{q} \le \min((n-1)(\frac{1}{2}-\frac{1}{r}),1) $ , $ n(\frac{1}{2}-\frac{1}{r})-\frac{1}{q} = \mu $ and $ A=-\Delta +1$. Then
\begin{equation}
\| e^{\pm itA^{1/2}} \phi_0 \|_{L^q_t(I,L^r_x({\bf R}^n))} \le c \| \phi_0 \|_{H^{\mu,2}_x({\bf R}^n)}
\label{-1.9}
\end{equation}
This implies
\begin{equation}
\label{-1.10}
\| f \|_{L^q_t(I,L^r_x({\bf R}^n))} \le c \| f\|_{X^{\mu,\frac{1}{2}+}_{\pm}(I)}
\end{equation}
and its dual version
\begin{equation}
\label{-1.11}
\| f \|_{X^{-\mu,-\frac{1}{2}-}_{\pm}(I)} \le c \| f \|_{L^{q'}_t(I,L^{r'}_x({\bf R}^n))}
\end{equation}
(see \cite{GTV}, Lemma 2.4 and \cite{GV}).\\
Finally we use the following notation for $\lambda \in {\bf R}: \, <\lambda>:=(1+\lambda^2)^{1/2} $, and $a+$ (resp. $a-$) denotes a number slightly larger (resp. smaller) than $a$.\\[2ex]
{\bf Acknowledgment:} I thank A. Gr\"unrock for many helpful discussions.  
\section{Bilinear estimates}
\begin{lemma}
\label{Lemma 0.0}
Assume $ P \in C^{\infty}({\bf R}^n,{\bf R}) \, , \, \psi \in C^{\infty}_0({\bf R}^n) $ with $ \nabla P \neq 0 $ on $ supp \, \psi $. Then the following identity holds:
$$ \int_{{\bf R}^n} \delta(P(x)) \psi(x) \, dx = \int_{\{P(x)=0\}} \frac{\psi(x)}{|\nabla P(x)|} \, dS_x $$
\end{lemma}
{\bf Proof:} We assume w.l.o.g. that on $ supp \, \psi $ we have $ \frac{\partial P}{\partial x_j} \neq 0 $ for some fixed $j$. Otherwise we take a smooth resolution of unity $ \{w_j\}^N_{j=1} $ s. th. on $ supp \, w_j $ we have this property. By \cite{GS}, p. 215 we have (where $H$ = Heaviside function):
\begin{eqnarray*}
\int_{{\bf R}^n} \delta(P(x))\psi(x)\,dx & = & \int_{{\bf R}^n} H'(P(x))\psi(x)\,dx \\
= \int_{{\bf R}^n} H'(P(x)) \frac{\partial P}{\partial x_j}(x) \frac{\psi(x)}{\frac{\partial P}{\partial x_j}(x)} \, dx & = & \int_{{\bf R}^n} \frac{\partial}{\partial x_j}(H(P(x)))  \frac{\psi(x)}{\frac{\partial P}{\partial x_j}(x)} \, dx \\
= - \int_{{\bf R}^n} H(P(x)) \frac{\partial}{\partial x_j} \left( \frac{\psi(x)}{\frac{\partial P}{\partial x_j}(x)} \right) \, dx & = & - \int_{\{P(x)>0\}}  \frac{\partial}{\partial x_j} \left( \frac{\psi(x)}{\frac{\partial P}{\partial x_j}(x)} \right) \, dx \\
= - \int_{\{P(x)=0\}} \frac{\psi(x)}{\frac{\partial P}{\partial x_j}(x)} \nu_j(x) \, dS_x & = & \int_{\{P(x)=0\}} \frac{\psi(x)}{|\nabla P(x)|} \, dS_x
\end{eqnarray*}
by Gauss' theorem, where $ \nu(x) = - \frac{\nabla P(x)}{|\nabla P(x)|} $ denotes the outer normal.\\[2ex]
Now, let $ \{\phi_j\}_{j=0}^{\infty} $ be a smooth resolution of unity in $ {\bf R}^3 $ , i.e. $ supp \, \phi_0 \subset \{| \xi|\le 2\} $, $ supp \, \phi_j \subset \{2^{j-1} \le | \xi| \le 2^{j+1} \} \, (j \ge 1) $ , $ \sum_{j=0}^{\infty} \phi_j (\xi) = 1 \, \forall \xi \in {\bf R}^3 $ , $ 0 \le \phi_j( \xi) \le 1 \, \forall j,\xi $ , $ \phi_j \in C^{\infty}({\bf R}^3).$
Define the operators $ P_{\Delta j} := {\cal F}^{-1} \phi_j {\cal F} \, , \, j=0,1,2,...\, .$ and $ P_l := \sum_{j=0}^l P_{\Delta j}$.\\
Let $A$ denote $-\Delta +1$ in the sequel.\\
The following generalized Strichartz inequality holds:
\begin{lemma}
\label{Lemma 0.1} For $l,m \in {\bf N}\cup \{0\}$:
$$ \| e^{it\Delta} P_{\Delta l} \psi_1 e^{\pm itA^{1/2}} P_{\Delta m}\psi_2 \|_{L^2_t({\bf R},L^2_x({\bf R}^3))} \le c 2^{m-\frac{l}{2}} \| P_{\Delta l} \psi_1 \|_{L^2_x({\bf R}^3)} \| P_{\Delta m} \psi_2 \|_{L^2_x({\bf R}^3)} $$
\end{lemma}
{\bf Proof:} Case 1: $ l \ge 2 $
\begin{eqnarray*}
&& \| e^{it\Delta} P_{\Delta l} \psi_1 e^{\pm itA^{1/2}} P_{\Delta m} \psi_2 \|_{L^2_{xt}}^2  \\ &
= & \int | \int_{\xi_1 + \xi_2 = \xi} e^{-it(|\xi_1|^2\pm (|\xi_2|^2+1)^{1/2})}  \phi_l(\xi_1) \widehat{\psi}_1(\xi_1)\phi_m(\xi_2)\widehat{\psi}_2(\xi_2) \, d\xi_1|^2d\xi dt  \\ & = &
\int_{\xi_1+\xi_2=\xi=\eta_1+\eta_2} e^{-it(|\xi_1|^2\pm (|\xi_2|^2+1)^{1/2}-|\eta_1|^2\mp (|\eta_2|^2+1)^{1/2})} \\ & & \phi_l(\xi_1)\widehat{\psi}_1(\xi_1) \phi_m(\xi_2)\widehat{\psi}_2(\xi_2)\phi_l(\eta_1)\overline{\widehat{\psi}}_1(\eta_1) \phi_m(\eta_2)\overline{\widehat{\psi}}_2(\eta_2)d\xi dt d\xi_1 d\eta_1 \\
& = & \int_{\xi_1 + \xi_2 = \xi = \eta_1 + \eta_2} \hspace{-2cm} \delta(|\xi_1|^2\pm (|\xi_2|^2+1)^{1/2}-|\eta_1|^2\mp (|\eta_2|^2+1)^{1/2})  \phi_l(\xi_1)\widehat{\psi}_1(\xi_1) \phi_m(\xi_2) \widehat{\psi}_2(\xi_2)\\
& & \widetilde{\phi}_l(\eta_1)^{1/2}\widetilde{\phi}_m(\eta_2)^{1/2}\phi_l(\eta_1)\overline{\widehat{\psi}}_1(\eta_1) \phi_m(\eta_2)\overline{\widehat{\psi}}_2(\eta_2)\widetilde{\phi}_l(\xi_1)^{1/2}\widetilde{\phi}_m(\xi_1)^{1/2}d\xi d\xi_1 d\eta_1 \\
& \le & \int_{\xi_1 + \xi_2 = \xi} |\phi_l(\xi_1)\widehat{\psi}_1(\xi_1)\phi_m(\xi_2)\widehat{\psi}_2(\xi_2)|^2 \\ & & ( \int_{\eta_1 + \eta_2 = \xi}\hspace{-1cm} \delta(|\eta_1|^2\pm (|\eta_2|^2+1)^{1/2}-|\xi_1|^2\mp (|\xi_2|^2+1)^{1/2}) \widetilde{\phi}_l(\eta_1) \widetilde{\phi}_m(\eta_2) d\eta_1) d\xi d\xi_1
\end{eqnarray*}
where $\widetilde{\phi}_l:=\phi_{l-1}+\phi_l+\phi_{l+1}=1$ on $supp \, \phi_l$.\\
Define now
$$ P(\eta_1):= |\eta_1|^2\pm (|\eta_2|^2+1)^{1/2}-|\xi_1|^2\mp (|\xi_2|^2+1)^{1/2} \, , \,
 \eta_2 = \xi - \eta_1 $$
Thus
$$ |\nabla_{\eta_1} P(\eta_1)| = | 2\eta_1 \pm \frac{\eta_1 - \xi}{(|\eta_1-\xi|^2+1)^{1/2}}| \ge 2|\eta_1| - 1 $$
Because $ supp \, \widetilde{\phi}_l \subset \{2^{l-2} \le |\eta_1| \le 2^{l+2} \} $ we have for $ l \ge 2 $:
$$ |\nabla_{\eta_1} P(\eta_1)| \ge 2^{l-1} - 1 \ge 2^{l-2} $$
Now using Lemma \ref{Lemma 0.0} we get
\begin{eqnarray*}
 I(\xi,\xi_1) & = & \int_{\eta_1 + \eta_2 = \xi} \delta(P(\eta_1)) \widetilde{\phi}_l(\eta_1) \widetilde{\phi}_m(\eta_2) d\eta_1 \\ 
& = & \int_ {P(\eta_1)=0 \, , \, \eta_1 + \eta_2 = \xi} \frac{\widetilde{\phi}_l(\eta_1) \widetilde{\phi}_m(\eta_2)}{|\nabla_{\eta_1} P(\eta_1)|} dS_{\eta_1} \\
& \le & 2^{-l+2} \int_{P(\eta_1) = 0} \widetilde{\phi}_m(\xi-\eta_1) \, dS_{\eta_1} \\
& \le & 2^{-l+2} \int_{\{P(\eta_1)=0\} \cap \{\eta_1 \in {\bf R}^3 \, : \,  |\xi-\eta_1| \le 2^{m+2}\}} dS_{\eta_1} \\
& \le & c  2^{-l+1}\, 2^{2(m+2)} \\
& = & c \, 2^{2m-l}
\end{eqnarray*} 
Thus
\begin{eqnarray*}
\lefteqn{\|e^{it\Delta} P_{\Delta l} \psi_1 e^{\pm itA^{1/2}} P_{\Delta m} \psi_2 \|_{L^2_{xt}}^2} \\
& & \le c \, 2^{2m-l} \int_{\xi_1 + \xi_2 = \xi} | \phi_l(\xi_1) \widehat{\psi}_1(\xi_1) \phi_m(\xi_2) \widehat{\psi}_2(\xi_2)|^2 \, d\xi d\xi_1 \\
& & = c \, 2^{2m-l} \int |\phi_l(\xi_1) \widehat{\psi_1}(\xi_1)|^2 \left(\int |\phi_m(\xi - \xi_1) \widehat{\psi}_2(\xi-\xi_1)|^2 \, d\xi \right) d\xi_1 \\
& & = c \, 2^{2m-l} \int | \widehat{P_{\Delta l} \psi_1}(\xi_1)|^2 \left( \int |(\widehat{P_{\Delta m} \psi_2})(\xi - \xi_1)|^2 d\xi \right) d\xi_1 \\
& & = c \, 2^{2m-l} \| P_{\Delta l} \psi_1 \|_{L^2_x}^2 \| P_{\Delta m} \psi_2 \|_{L^2_x}^2 
\end{eqnarray*}
Case 2: $ l\le 1 $ \\
\begin{eqnarray*}
\|e^{it\Delta} P_{\Delta l} \psi_1 e^{\pm itA^{1/2}} P_{\Delta m} \psi_2 \|_{L^2_{xt}} 
& \le & \| e^{it\Delta} P_{\Delta l} \psi_1\|_{L^4_{xt}} \| e^{\pm itA^{1/2}} P_{\Delta m} \psi_2 \|_{L^4_{xt}} \\
& \le & c \| P_{\Delta l} \psi_1 \|_{\dot{H}^{\frac{1}{4},2}_x({\bf R}^3)} \| P_{\Delta m} \psi_2 \|_{H^{\frac{1}{2},2}_x({\bf R}^3)}
\end{eqnarray*}
using Strichartz` inequalities as follows: by (\ref{-1.6}) with $ q=4$ , $ r=3$ we have
$$ \| e^{it\Delta} P_{\Delta l} \psi_1\|_{L^4_{xt}} \le c \|e^{it\Delta} P_{\Delta l} \psi_1\|_{L^4_t\dot{H}^{\frac{1}{4},3}_x} \le c \|P_{\Delta l} \psi_1\|_{\dot{H}^{{\frac{1}{4}},2}_x({\bf R}^3)} $$
and by (\ref{-1.9}) with $ q=r=4$ , $ \mu =-1/2$:
$$ \|e^{\pm itA^{1/2}} P_{\Delta m} \psi_2 \|_{L^4_t L^4_x} \le c \| P_{\Delta m} \psi_2 \|_{H^{\frac{1}{2},2}_x({\bf R}^3)} $$
Now $ \|P_{\Delta l} \psi_1 \|_{\dot{H}^{\frac{1}{4},2}} \le c \|P_{\Delta l} \psi_1 \|_{L^2_x} $ and $ \|P_{\Delta m} \psi_2 \|_{H^{\frac{1}{2},2}} \le c 2^{\frac{m}{2}} \|P_{\Delta m} \psi_2 \|_{L^2_x} $ gives the claimed result.
\begin{lemma}
\label{Lemma 0.2}
If $ 0 \le s < 1/2 $ the following estimate holds:
$$ \| e^{it\Delta} u_1 e^{\pm itA^{1/2}} u_2 \|_{L^2_t({\bf R},H^{s,2}_x({\bf R}^3))} \le c \|u_1\|_{L^2_x({\bf R}^3)} \|u_2\|_{H^{s+\frac{1}{2}+,2}({\bf R}^3)} $$
\end{lemma}
{\bf Proof:}
With $ P_{-1} := 0 $ we have
\begin{eqnarray*}
\lefteqn{\|e^{it\Delta} u_1 e^{\pm itA^{1/2}}u_2\|_{L^2_tH^{s,2}_x} = \lim_{N \to \infty} \| e^{it\Delta} P_N u_1 e^{\pm itA^{1/2}} P_N u_2 \|_{L^2_t H^{s,2}_x}} \\
& & \le \lim_{N \to \infty} \sum_{l=0}^N \|e^{it\Delta}P_l u_1 e^{\pm itA^{1/2}} P_l u_2 - e^{it\Delta}P_{l-1}u_1 e^{\pm itA^{1/2}} P_{l-1}u_2 \|_{L^2_t H^{s,2}_x} \\
& & \le \sum_{l=0}^{\infty} \|e^{it\Delta} P_{\Delta l} u_1 e^{\pm itA^{1/2}}P_lu_2\|_{L^2_t H^{s,2}_x} + \sum_{l=0}^{\infty} \|e^{it\Delta} P_{l-1} u_1 e^{\pm itA^{1/2}}P_{\Delta l}u_2\|_{L^2_t H^{s,2}_x} \\
& & =: \Sigma_1 + \Sigma_2
\end{eqnarray*}
Now 
\begin{eqnarray*}
\lefteqn{{\cal F}(e^{it\Delta} P_{\Delta l} u_1 e^{\pm itA^{1/2}}P_lu_2)(\xi) = [{\cal F}(e^{it\Delta} P_{\Delta l} u_1) \ast {\cal F}(e^{\pm itA^{1/2}} P_l u_2)](\xi) } \\
& & = \int \int_{\xi_1 + \xi_2 = \xi} {\cal F}(e^{it\Delta} P_{\Delta l} u_1)(\xi_1)  {\cal F}(e^{\pm itA^{1/2}} P_l u_2)(\xi_2)\, d\xi_1 d\xi_2
\end{eqnarray*}
where $ |\xi| \le 2^{l+2} $ because of $ 2^{l-1} \le |\xi_1| \le 2^{l+1} $ and $|\xi_2| \le 2^{l+1} $ so that \\ $ supp \, {\cal F}(e^{it\Delta}P_{\Delta l} u_1 e^{\pm itA^{1/2}} P_l u_2) \subset \{ |\xi| \le 2^{l+2} \} $.\\
Thus
\begin{eqnarray*}
\|e^{it\Delta} P_{\Delta l} u_1 e^{\pm itA^{1/2}} P_l u_2 \|_{L^2_tH^{s,2}_x} & = & \| <\xi>^s {\cal F}(e^{it\Delta} P_{\Delta l} u_1 e^{\pm itA^{1/2}} P_l u_2) \|_{L^2_t L^2_{\xi}} \\
& \le & c 2^{sl} \|e^{it\Delta} P_{\Delta l} u_1 e^{\pm itA^{1/2}} P_l u_2\|_{L^2_{xt}}
\end{eqnarray*}
so that by Lemma \ref{Lemma 0.1}:
\begin{eqnarray*}
\Sigma_1 & \le & c \sum_{l=0}^{\infty} 2^{sl} \|e^{it\Delta} P_{\Delta l} u_1 e^{\pm itA^{1/2}} P_l u_2\|_{L^2_{xt}} \\
& \le & c \sum_{l=0}^{\infty} 2^{sl} \sum_{n \le l} \|e^{it\Delta} P_{\Delta l} u_1 e^{\pm itA^{1/2}} P_{\Delta n} u_2\|_{L^2_{xt}} \\
& \le & c \sum_{l=0}^{\infty} 2^{sl} \sum_{n \le l} 2^{n-\frac{l}{2}} \| P_{\Delta l} u_1 \|_{L^2_x} \| P_{\Delta n} u_2\|_{L^2_x} \\
& \le & c \sum_{l=0}^{\infty} 2^{l(s-\frac{1}{2})} \sum_{n \le l} 2^n \| u_1 \|_{L^2_x} \| P_{\Delta n} u_2\|_{L^2_x} \\
& = & c \sum_{l=0}^{\infty} 2^{l(s-\frac{1}{2})} \sum_{n \le l} 2^n 2^{-n(s+\frac{1}{2}+)} 2^{n(s+\frac{1}{2}+)}\| u_1 \|_{L^2_x} \| P_{\Delta n} u_2\|_{L^2_x} \\
& \le & c \sum_{l=0}^{\infty} 2^{l(s-\frac{1}{2})} \sum_{n \le l} 2^{n(\frac{1}{2}-s-)}\| u_1 \|_{L^2_x} \| u_2\|_{H^{s+\frac{1}{2}+,2}_x} \\
& \le & c \sum_{l=0}^{\infty} 2^{l(s-\frac{1}{2})} 2^{l(\frac{1}{2}-s-)}\| u_1 \|_{L^2_x} \| u_2\|_{H^{s+\frac{1}{2}+,2}_x} \quad (\mbox{if} \, s < \frac{1}{2}) \\
& = & c \sum_{l=0}^{\infty} (2^{0-})^l \| u_1 \|_{L^2_x} \| u_2\|_{H^{s+\frac{1}{2}+,2}_x} \\
& \le & c \| u_1 \|_{L^2_x} \| u_2\|_{H^{s+\frac{1}{2}+,2}_x}
\end{eqnarray*}
In order to estimate $ \Sigma_2 $ we use Strichartz' inequalities again to conclude as above:
\begin{eqnarray*}
\Sigma_2 
& \le & c \sum_{l=0}^{\infty} 2^{sl} \| e^{it\Delta} P_{l-1} u_1 e^{\pm itA^{1/2}} P_{\Delta l} u_2 \|_{L^2_{xt}} \\
& \le & c \sum_{l=0}^{\infty} 2^{sl} \| e^{it\Delta} P_{l-1} u_1\|_{L^{2+}_t L^{6-}_x} \| e^{\pm itA^{1/2}} P_{\Delta l} u_2 \|_{L^{\infty -}_t L^{3+}_x} \\
& \le & c \sum_{l=0}^{\infty} 2^{sl} \|  P_{l-1} u_1\|_{L^2_x} \|  P_{\Delta l} u_2 \|_{H^{\frac{1}{2}+,2}_x} \\
& \le & c \sum_{l=0}^{\infty} (2^{0-})^l 2^{l(s+)} \| u_1\|_{L^2_x} \|  P_{\Delta l} u_2 \|_{H^{\frac{1}{2}+,2}_x} \\
& \le & c \sum_{l=0}^{\infty} (2^{0-})^l \| u_1\|_{L^2_x} \| u_2 \|_{H^{s+\frac{1}{2}+,2}_x} \\
& \le & c \| u_1\|_{L^2_x} \| u_2 \|_{H^{s+\frac{1}{2}+,2}_x}
\end{eqnarray*}
This completes the proof.

The following lemma is a direct consequence and a bilinear version of \cite{GTV}, Lemma 2.3.
\begin{lemma}
\label{Lemma 0.3}
For $0 \le s < 1/2 $ the following estimate holds:
$$ \| v_1 v_2 \|_{L^2_t({\bf R},H^{s,2}_x({\bf R}^3))} \le c \|v_1\|_{X^{0,\frac{1}{2}+}} \|v_2\|_{X_{\pm}^{s+\frac{1}{2}+,\frac{1}{2}+}} $$
\end{lemma}
{\bf Proof:}
We define $ U_{\varphi_1}(t):= e^{it\Delta} \, , \, U_{\varphi_2}(t):= e^{\pm itA^{1/2}} $ and start from
$$ v_j(t) = U_{\varphi_j}(t) \int e^{it\tau} ({\cal F}_t U_{\varphi_j}(-t)v_1)(\tau)\, d\tau $$
Thus by use of Lemma \ref{Lemma 0.2} we get
\begin{eqnarray*}
\lefteqn{ \|v_1 v_2\|_{L^2_t({\bf R},H^{s,2}_x)}} \\
& & \hspace{-0.9cm} = \| \int \int e^{it\tau_1} U_{\varphi_1}(t)({\cal F}_t U_{\varphi_1}(-t)v_1)(\tau_1) e^{it\tau_2} U_{\varphi_2}(t)({\cal F}_t U_{\varphi_2}(-t)v_2)(\tau_2) \, d\tau_1 d\tau_2 \|_{L^2_t H^{s,2}_x} \\
& & \hspace{-0.9cm} \le \int \int \| U_{\varphi_1}(t)({\cal F}_t U_{\varphi_1}(-t)v_1)(\tau_1) U_{\varphi_2}(t)({\cal F}_t U_{\varphi_2}(-t)v_2)(\tau_2)  \|_{L^2_t H^{s,2}_x}d\tau_1d\tau_2 \\
& & \hspace{-0.9cm} \le c \int \int \| ({\cal F}_t U_{\varphi_1}(-t)v_1)(\tau_1)\|_{L^2_x} \|({\cal F}_t U_{\varphi_2}(-t)v_2)(\tau_2)  \|_{ H^{s+\frac{1}{2}+,2}_x}d\tau_1d\tau_2 \\
& & \hspace{-0.9cm} \le c (\int <\tau_1>^{-1-} d\tau_1)^{\frac{1}{2}}(\int <\tau_1>^{1+} \|({\cal F}_t U_{\varphi_1}(-t)v_1)(\tau_1)\|_{L^2_x}^2 d\tau_1)^{\frac{1}{2}} \\
& & \hspace{-0.9cm} \quad \cdot (\int <\tau_2>^{-1-} d\tau_2)^{\frac{1}{2}}(\int <\tau_2>^{1+} \|({\cal F}_t U_{\varphi_2}(-t)v_2)(\tau_2)\|_{H^{s+\frac{1}{2}+,2}_x}^2 d\tau_2)^{\frac{1}{2}} \\
& & \hspace{-0.9cm} \le c \| v_1 \|_{X^{0,\frac{1}{2}+}} \| v_2 \|_{X_{\pm}^{s+\frac{1}{2}+,\frac{1}{2}+}}
\end{eqnarray*}
\begin{Cor}
\label{Corollary 0.4}
The following estimates hold for $ 0 \le s < 1/2 $:
$$ \| un\|_{X^{0,-\frac{1}{2}-}} \le c \|u\|_{L^2_t({\bf R},H^{-s,2}_x)} \|n\|_{X^{s+\frac{1}{2}+,\frac{1}{2}+}_{\pm}} $$
Here $u$ and/or $n$ can be replaced by $\overline{u}$ and/or $\overline{n}$ on the left and/or right hand side.
\end{Cor}
{\bf Proof:}
By Lemma \ref{Lemma 0.3} the mapping $ X^{0,\frac{1}{2}+} \to L^2_t H^{s,2}_x $ defined by $ u \mapsto un $ is bounded by $ c\|n\|_{X_{\pm}^{s+\frac{1}{2}+,\frac{1}{2}+}} $. Thus the dual mapping $ L^2_tH^{-s,2}_x \to X^{0,-\frac{1}{2}-} $ defined by $ u \mapsto u \overline{n} $ has the same bound, i.e.
$$ \|u\overline{n}\|_{X^{0,-\frac{1}{2}-}} \le c \|u\|_{L^2_t H^{-s,2}_x} \|n\|_{X^{s+\frac{1}{2}+,\frac{1}{2}+}_{\pm}} $$
Because
\begin{eqnarray*}
 \| \overline{n} \|_{X^{s+\frac{1}{2}+,\frac{1}{2}+}_+}^2 & 
 = & \int \int | <\tau + |\xi|>^{\frac{1}{2}+} <\xi>^{s+\frac{1}{2}+} \widehat{\overline{n}}(\xi,\tau)|^2 d\xi d\tau \\
& = &  \int \int <\tau + |\xi|>^{\frac{1}{2}+} <\xi>^{s+\frac{1}{2}+} \overline{\widehat{n}}(-\xi,-\tau)|^2 d\xi d\tau \\
& = & \int \int <\widetilde{\tau} - |\widetilde{\xi}|>^{\frac{1}{2}+} <\widetilde{\xi}>^{s+\frac{1}{2}+} \widehat{n}(\widetilde{\xi},\widetilde{\tau})|^2 d\widetilde{\xi} d\widetilde{\tau} \\
& = & \|n\|_{X^{s+\frac{1}{2}+,\frac{1}{2}+}_-}^2
\end{eqnarray*}
we can replace $n$ by $\overline{n}$ on the left and/or right hand side of the claimed inequality. Trivially $u$ may be replaced by $\overline{u}$ on the left and/or right hand side.\\[2ex]
The key estimate for the nonlinearity in the Schr\"odinger equation is given in the following
\begin{lemma}
\label{Lemma 0.5}
If $ s \ge 0 \, , \, 0 \le \sigma < 1/2 $ the following estimates hold:
$$ \|un_{\pm}\|_{X^{s,-\frac{1}{2}-}} \le c(\|u\|_{X^{0,\frac{1}{2}+}} \|n_{\pm}\|_{X_{\pm}^{s-\frac{1}{2},\frac{1}{2}+}} + \|u\|_{X^{s-\sigma,0}} \|n_{\pm}\|_{X_{\pm}^{\sigma+\frac{1}{2}+,\frac{1}{2}+}}) $$
Here $u$ and/or $n_{\pm}$ may be replaced by $\overline{u}$ and/or $\overline{n}_{\pm}$ on the left and/or right hand side.
\end{lemma}
{\bf Proof:}
We take the scalar product with a function $w \in X^{-s,\frac{1}{2}+}$ and estimate
$$ | \int (un_{\pm})(x_1,t_1)w(x_1,t_1)dx_1dt_1| = |\int \int \widehat{u}(\xi_2,\tau_2)\widehat{n_{\pm}}(\xi,\tau)\widehat{w}(\xi_1,\tau_1)d\xi_1d\tau_1d\xi_2d\tau_2| $$
where $ \xi = \xi_1 - \xi_2 \, , \, \tau = \tau_1 - \tau_2$. \\
We split the integral into the parts $ B_1 \, : \, |\xi_1| \le 2 |\xi_2| $ and $ B_2 \, : \, |\xi_1| > 2 |\xi_2| $.\\
Estimate of $B_1$:\\
We use the notation $ \sigma_i = \tau_i + |\xi_i|^2 \, (i=1,2) $, $ \sigma_{\pm}=\tau \pm |\xi| $ and Corollary \ref{Corollary 0.4}: 
\begin{eqnarray*}
\lefteqn{|\int_{B_1}\int \widehat{u}(\xi_2,\tau_2)\widehat{n_{\pm}}(\xi,\tau)\widehat{w}(\xi_1,\tau_1)d\xi_1d\tau_1d\xi_2d\tau_2|} \\ 
& & = |\int <\xi_1>^s <\sigma_1>^{-\frac{1}{2}-}(\int_{|\xi_1|\le 2|\xi_2|} \widehat{u}(\xi_2,\tau_2)\widehat{n_{\pm}}(\xi,\tau)d\xi_2d\tau_2) \cdot  \\*[0cm] 
& & \qquad \qquad \qquad \qquad \qquad \qquad \qquad \quad  \cdot <\xi_1>^{-s}<\sigma_1>^{\frac{1}{2}+}\widehat{w}(\xi_1,\tau_1)d\xi_1d\tau_1|
 \\
& & \le \| <\xi_1>^s <\sigma_1>^{-\frac{1}{2}-}\int_{|\xi_1|\le 2|\xi_2|} \widehat{u}(\xi_2,\tau_2) \widehat{n_{\pm}}(\xi,\tau)d\xi_2d\tau_2\|_{L^2_{\xi_1\tau_1}}\|w\|_{X^{-s,\frac{1}{2}+}} \\
& & \le c \| <\sigma_1>^{-\frac{1}{2}-}\int_{|\xi_1|\le 2|\xi_2|} <\xi_2>^s \widehat{u}(\xi_2,\tau_2) \widehat{n_{\pm}}(\xi,\tau)d\xi_2d\tau_2\|_{L^2_{\xi_1\tau_1}}\|w\|_{X^{-s,\frac{1}{2}+}} \\
& & = c \| (A^{s/2} u) n_{\pm}\|_{X^{0,-\frac{1}{2}-}} \|w\|_{X^{-s,\frac{1}{2}+}} \\
& & \le c \|A^{s/2} u\|_{L_t^2 H_x^{-\sigma,2}} \| n_{\pm}\|_{X_{\pm}^{\sigma+\frac{1}{2}+,\frac{1}{2}+}} \|w\|_{X^{-s,\frac{1}{2}+}} \\
& & \le c \| u\|_{X^{s-\sigma,0}} \| n_{\pm}\|_{X_{\pm}^{\sigma+\frac{1}{2}+,\frac{1}{2}+}} \|w\|_{X^{-s,\frac{1}{2}+}} 
\end{eqnarray*}
Estimate of $B_2$: \\
Define
$$ \widehat{v_2}:=<\sigma_2>^{\frac{1}{2}+}\widehat{u} \, , \, \widehat{v}:=<\xi>^{s-\frac{1}{2}}<\sigma_{\pm}>^{\frac{1}{2}+}\widehat{n_{\pm}} \, , \, \widehat{v_1}:=<\xi_1>^{-s}<\sigma_1>^{\frac{1}{2}+}\widehat{w} $$
so that
$$ \|u\|_{X^{0,\frac{1}{2}+}} = \|v_2\|_{L^2_{xt}} \, , \, \|n_{\pm}\|_{X_{\pm}^{s-\frac{1}{2},\frac{1}{2}+}} = \|v\|_{L^2_{xt}} \, , \, \|w\|_{X^{-s,\frac{1}{2}+}} = \| v_1\|_{L^2_{xt}} $$
In $B_2$ we have
$$ \frac{1}{2}|\xi_1| \le |\xi_1| - |\xi_2| \le |\xi| \le |\xi_1|+|\xi_2| \le \frac{3}{2} |\xi_1| $$
and
\begin{eqnarray*}
z_{\pm} & := & |\xi_1|^2 - |\xi_2|^2 \mp |\xi| = (\sigma_1 - \tau_1)-(\sigma_2 - \tau_2) + (\tau - \sigma_{\pm}) \\
& = & \sigma_1 - \sigma_2 - \sigma_{\pm} + \tau_2 - \tau_1 + \tau = \sigma_1 - \sigma_2 - \sigma_{\pm}
\end{eqnarray*}
Thus
\begin{eqnarray*}
\left| |\xi_1|^2 - |\xi_2|^2 \mp |\xi| \right| & \ge & |\xi_1|^2 - |\xi_2|^2 - |\xi| \ge |\xi_1|^2 - \frac{1}{4}|\xi_1|^2 - \frac{3}{2} |\xi_1| \\
& \ge & \frac{3}{4} |\xi_1|^2 - \frac{3}{8}|\xi_1|^2 - 10 = \frac{3}{8} |\xi_1|^2 - 10
\end{eqnarray*}
This implies
$$ \frac{3}{8} |\xi_1|^2 \le \left| |\xi_1|^2 - |\xi_2|^2 \mp |\xi| \right| + 10 \le |\sigma_1| + |\sigma_2| + |\sigma_{\pm}| + 10 $$ 
and
$$ <\xi_1>^{\frac{1}{2}} \le c \left( <\sigma_1>^{\frac{1}{4}} + <\sigma_2>^{\frac{1}{4}} + <\sigma_{\pm}>^{\frac{1}{4}} \right) $$
Therefore we get
\begin{eqnarray*}
\lefteqn{ \left| \int_{B_2}\int \widehat{u}(\xi_2,\tau_2) \widehat{n_{\pm}}(\xi,\tau) \widehat{w}(\xi_1,\tau_1)d\xi_1 d\tau_1 d\xi_2 d\tau_2 \right| } \\
& & = \left| \int_{B_2} \int \frac{\widehat{v_2} \widehat{v} \widehat{v_1} <\xi_1>^s}{<\sigma_2>^{\frac{1}{2}+} <\sigma_{\pm}>^{\frac{1}{2}+} <\xi>^{s-\frac{1}{2}} <\sigma_1>^{\frac{1}{2}+}} d\xi_1 d\tau_1 d\xi_2 d\tau_2 \right| \\
& & \le \int \int \frac{| \widehat{v_2} \widehat{v} \widehat{v_1}| <\xi_1>^{\frac{1}{2}}}{<\sigma_2>^{\frac{1}{2}+} <\sigma_{\pm}>^{\frac{1}{2}+}  <\sigma_1>^{\frac{1}{2}+}} d\xi_1 d\tau_1 d\xi_2 d\tau_2 \\
& & \le \int \int \frac{| \widehat{v_2} \widehat{v} \widehat{v_1}| \left( <\sigma_1>^{\frac{1}{4}} + <\sigma_2>^{\frac{1}{4}} + <\sigma_{\pm}>^{\frac{1}{4}} \right)}{<\sigma_2>^{\frac{1}{2}+} <\sigma_{\pm}>^{\frac{1}{2}+}  <\sigma_1>^{\frac{1}{2}+}} 
d\xi_1 d\tau_1 d\xi_2 d\tau_2 \\
& & \le I_1 + I_2 + I_3
\end{eqnarray*}
where
\begin{eqnarray*}
I_1 & = & \int \int \frac{| \widehat{v} \widehat{v_1} \widehat{v_2}| }{<\sigma_{\pm}>^{\frac{1}{4}+} <\sigma_1>^{\frac{1}{2}+}  <\sigma_2>^{\frac{1}{2}+}} 
d\xi_1 d\tau_1 d\xi_2 d\tau_2 \\
I_2 & = & \int \int \frac{| \widehat{v} \widehat{v_1} \widehat{v_2}| }{<\sigma_{\pm}>^{\frac{1}{2}+} <\sigma_1>^{\frac{1}{4}+}  <\sigma_2>^{\frac{1}{2}+}} 
d\xi_1 d\tau_1 d\xi_2 d\tau_2 \\
I_3 & = & \int \int \frac{| \widehat{v} \widehat{v_1} \widehat{v_2}| }{<\sigma_{\pm}>^{\frac{1}{2}+} <\sigma_1>^{\frac{1}{2}+}  <\sigma_2>^{\frac{1}{4}+}} 
d\xi_1 d\tau_1 d\xi_2 d\tau_2 
\end{eqnarray*}
Our aim is to estimate each of these integrals by
$$ c \|v\|_{L^2_{xt}}\|v_1\|_{L^2_{xt}}\|v_2\|_{L^2_{xt}} = c \|n_{\pm}\|_{X_{\pm}^{s-\frac{1}{2},\frac{1}{2}+}} \|u\|_{X^{0,\frac{1}{2}+}}\|w\|_{X^{-s,\frac{1}{2}+}} $$
These estimates together with the estimate of $B_1$ imply the claimed results.\\
$I_2$ and $I_3$ can be treated in exactly the same way, so we are left with $I_1$ and $I_2$.\\
Estimate of $I_1$:\\
By H\"older's inequality we get
$$ I_1 \le \| {\cal F}^{-1}(<\sigma_{\pm}>^{-\frac{1}{4}-} |\widehat{v}|)\|_{L^4_t L^2_x} \prod_{i=1}^2 \|{\cal F}^{-1}(<\sigma_i>^{-\frac{1}{2}-}|\widehat{v_i}|)\|_{L^{8/3}_t L^4_x} $$
Using Strichartz' inequality for the Schr\"odinger equation we get from (\ref{-1.7}) with $ f := {\cal F}^{-1}(<\sigma_i>^{-\frac{1}{2}-}|\widehat{v_i}|) $: 
$$ \|{\cal F}^{-1}(<\hspace{-0.1cm} \sigma_i\hspace{-0.1cm}>^{-\frac{1}{2}-}|\widehat{v_i}|)\|_{L^{8/3}_t L^4_x} \le c \|{\cal F}^{-1}(<\hspace{-0.1cm}\sigma_i\hspace{-0.1cm}>^{-\frac{1}{2}-}|\widehat{v_i}|)\|_{X^{0,\frac{1}{2}+}} = \|v_i\|_{X^{0,0}} = \|v_i\|_{L^2_{xt}} $$
Similarly using the trivial estimate $ \|e^{\pm itA^{1/2}} u_0\|_{L^{\infty}_t L^2_x} \le c \|u_0\|_{L^2_x} $ we get \\ $ \|w\|_{L^{\infty}_t L^2_x} \le c \|w\|_{X_{\pm}^{0,\frac{1}{2}+}}$. Interpolating with $ \|w\|_{L^2_tL^2_x} = \|w\|_{X_{\pm}^{0,0}} $ we arrive at \\ $ \|w\|_{L^4_tL^2_x} \le  c \|w\|_{X_{\pm}^{0,\frac{1}{4}+}} $. With $ w:= {\cal F}^{-1}(<\sigma_{\pm}>^{-\frac{1}{4}-}|\widehat{v}|) $ this implies
$$ \|{\cal F}^{-1}(<\sigma_{\pm}>^{-\frac{1}{4}-}|\widehat{v}|)\|_{L^4_tL^2_x} \le c \| {\cal F}^{-1}(<\sigma_{\pm}>^{-\frac{1}{4}-}|\widehat{v}|)\|_{X_{\pm}^{0,\frac{1}{4}+}} = \|v\|_{X_{\pm}^{0,0}} = \|v\|_{L^2_{xt}} $$
Thus $I_1$ is estimated in the desired way.\\
Estimate of $I_2$:\\
We have
\begin{eqnarray*}
I_2 & \le & \|{\cal F}^{-1}(<\sigma_{\pm}>^{-\frac{1}{2}-\epsilon}|\widehat{v}|)\|_{L^{\infty -}_t L^2_x} \|{\cal F}^{-1}(<\sigma_1>^{-\frac{1}{4}-\epsilon}|\widehat{v_1}|)\|_{L^{2+}_t L^{3+}_x} \cdot \\
& &
\cdot \|{\cal F}^{-1}(<\sigma_2>^{-\frac{1}{2}-\epsilon}|\widehat{v_2}|)\|_{L^{2+}_t L^{6-}_x}
\end{eqnarray*}
Now, as above, $\|w\|_{L^{\infty}_t L^2_x} \le c \|w\|_{X_{\pm}^{0,\frac{1}{2}+}}$. Interpolating with $\|w\|_{L^2_t L^2_x} = \|w\|_{X^{0,0}_{\pm}}$ we get $\|w\|_{L^{\infty -}_t L^2_x} \le c \|w\|_{X^{0,\frac{1}{2}+}_{\pm}}$. With $w:={\cal F}^{-1}(<\sigma_{\pm}>^{-\frac{1}{2}-\epsilon}|\widehat{v}|)$ this implies
$$ \|{\cal F}^{-1}(<\hspace{-0.15cm}\sigma_{\pm}\hspace{-0.15cm}>^{-\frac{1}{2}-\epsilon}|\widehat{v}|)\|_{L^{\infty -}_t L^2_x} \le c \|{\cal F}^{-1}(<\hspace{-0.15cm}\sigma_{\pm}\hspace{-0.15cm}>^{-\frac{1}{2}-\epsilon}|\widehat{v}|)\|_{X^{0,\frac{1}{2}+}_{\pm}} \le c \|v\|_{X^{0,0}_{\pm}} = c \|v\|_{L^2_{xt}} $$ 
Using Strichartz' estimate for the Schr\"odinger equation we get from (\ref{-1.7}):
$$ \|{\cal F}^{-1}(<\sigma_2>^{-\frac{1}{2}-\epsilon}|\widehat{v_2}|)\|_{L^{2+}_t L^{6-}_x} \le c \|{\cal F}^{-1}(<\sigma_2>^{-\frac{1}{2}-\epsilon}|\widehat{v_2}|)\|_{X^{0,\frac{1}{2}+}} \le c \|v_2\|_{L^2_{xt}} $$
and, finally, from $ \|f\|_{L^{2+}_t L^{6-}_x} \le c \|f\|_{X^{0,\frac{1}{2}+}}$ and the trivial identity $\|f\|_{L^2_t L^2_x} = \|f\|_{X^{0,0}}$ we get by interpolation $\|f\|_{L^{2+}_t L^{3+}_x} \le c \|f\|_{X^{0,\frac{1}{4}+}} $. \\With $ f = {\cal F}^{-1}(<\sigma_1>^{-\frac{1}{4}-\epsilon}|\widehat{v_1}|) $ this implies
\begin{eqnarray*}
\|{\cal F}^{-1}(<\sigma_1>^{-\frac{1}{4}-\epsilon}|\widehat{v_1}|)\|_{L^{2+}_t L^{3+}_x} & \le & c \|{\cal F}^{-1}(<\sigma_1>^{-\frac{1}{4}-\epsilon}|\widehat{v_1}|)\|_{X^{0,\frac{1}{4}+}} \\
& \le & c \|v_1\|_{X^{0,0}} = c \|v_1\|_{L^2_{xt}}
\end{eqnarray*}
which completes the proof.
\section{Local existence and uniqueness}
The system (\ref{-1.1}),(\ref{-1.2}),(\ref{-1.3}) is transformed into a system of first order in $t$ in the usual way. Defining $ A:= -\Delta +1 $ and 
$$ \phi_{\pm} := \phi \pm iA^{-\frac{1}{2}}\phi_t $$
we have
\begin{eqnarray}
\label{1.3a}
\phi & = & \frac{1}{2} (\phi_+ + \phi_-) \\
\label{1.3b}
2iA^{-\frac{1}{2}}\phi_t & = & \phi_+ - \phi_-
\end{eqnarray}
and the equivalent system is
\begin{eqnarray}
\label{1.4}
i \psi_t + \Delta \psi & = & -\frac{1}{2}(\phi_+ + \phi_-)\psi \\
\label{1.5}
i \phi_{\pm t} \mp A^{\frac{1}{2}}\phi_{\pm} & = & \mp A^{-\frac{1}{2}}(|\psi|^2) \\
\psi(0) = \psi_0 \, , \, \phi_{\pm}(0) & = & \phi_0 \pm iA^{-\frac{1}{2}}\phi_1 = : \phi_{0\pm}
\end{eqnarray}
The corresponding system of integral equations reads as follows:
\begin{eqnarray}
\label{1.6}
\psi(t) & = & e^{it\Delta}\psi_0 + i \int^t_0 e^{i(t-s)\Delta} \frac{1}{2}(\phi_+(s) + \phi_-(s)) \psi(s) \, ds \\
\label{1.7}
\phi_{\pm}(t) & = & e^{\mp itA^{\frac{1}{2}}}\phi_{0\pm} \pm i \int^t_0 e^{\mp i(t-s)A^{\frac{1}{2}}} A^{-\frac{1}{2}}(|\psi(s)|^2) \, ds
\end{eqnarray}
We always assume $t \in I = [0,|I|]$. In this case we could, whenever helpful, place a factor $\psi_1(t)$ in front of the first terms on the right hand side and a factor $\psi_{|I|}(t)$ in front of any of the integrals without changing the equations at all. Here $\psi \in C^{\infty}_0({\bf R})$ is a nonnegative cutoff function with $\psi(t)=0$ if $|t|\ge 2$, $\psi(t)=1$ if $|t|\le 1$ and $\psi_{\delta}(t):=\psi(\frac{t}{\delta})$.\\
The decisive estimates for the nonlinearities follow from the corresponding results for the Zakharov system in \cite{GTV} as follows.
\begin{lemma}
\label{Lemma 1.1}
The estimate
\begin{equation}
\label{1.8}
\| |\psi|^2 \|_{X^{m-1,-\frac{1}{2}+}_{\pm}} \le c \| \psi \|^2_{X^{s,\frac{1}{2}+}}
\end{equation}
 holds, provided $ s \ge 0 \, , \, 2s \ge m-1 \, , \, s > m-2 $.
\end{lemma}
{\bf Proof:}
follows from the proof of \cite{GTV}, Lemma 3.5 with $l=m-2$, $k=s$, $c=\frac{1}{2}-$, $b_1 = \frac{1}{2}+$.
\begin{lemma}
\label{Lemma 1.2}
Let $ m \ge 0 \, , \, s < m+1 $. Then the following estimates hold for $\epsilon > 0$:
\begin{eqnarray}
\label{1.9}
\|\phi_{\pm}\psi\|_{X^{s,-\frac{1}{2}+}} & \le & c \|\phi_{\pm}\|_{X^{m,\frac{1}{2}+}_{\pm}} \|\psi\|_{X^{s,\frac{1}{2}+}} \\
\label{1.10}
\|\phi_{\pm}\psi\|_{X^{0,-\frac{1}{4}-\epsilon}} & \le & c \|\phi_{\pm}\|_{X^{0,\frac{1}{2}+\epsilon}_{\pm}} \|\psi\|_{X^{0,\frac{1}{2}+\epsilon}}
\end{eqnarray}
\end{lemma}
{\bf Proof:} \cite{GTV}, Lemma 3.4.\\[2ex]
The following local wellposedness result is a direct consequence of (\ref{1.8}) and (\ref{1.9}).
\begin{theorem}
\label{Theorem 1.3} Let $s$ and $m$ satisfy $s,m\ge 0 \, , \, m-2<s<m+1 \, , \, 2s\ge m-1$. Then the system (\ref{-1.1}),(\ref{-1.2}),(\ref{-1.3}) with initial data $$ (\psi_0,\phi_0,\phi_1) \in H^{s,2}({\bf R}^3) \times H^{m,2}({\bf R}^3)
\times H^{m-1,2}({\bf R}^3) $$ has a unique local solution
\begin{eqnarray*}
\lefteqn{(u,\phi,\phi_t)} \\ & & \hspace{-0.8cm} \in X^{s,\frac{1}{2}+}[0,T] \times (X^{m,\frac{1}{2}+}_+[0,T] + X^{m,\frac{1}{2}+}_-[0,T]) \times (X^{m-1,\frac{1}{2}+}_+[0,T] + X^{m-1,\frac{1}{2}+}_-[0,T]) 
\end{eqnarray*}
We also have
$$ u \in C^0([0,T],H^{s,2}({\bf R}^3)) \, , \, \phi \in C^0([0,T],H^{m,2}({\bf R}^3))\cap C^1([0,T],H^{m-1,2}({\bf R}^3)) $$
\end{theorem}
{\bf Proof:} as in \cite{GTV}.\\
{\bf Remark:} Especially the range $0\le s \le 1$ , $ 0<m\le 1$ is admissible as well as $s=m=0$.
\section{Energy bounds and decomposition of data}
It is well known that the following conservation laws hold for the system (\ref{-1.1}),(\ref{-1.2}),\\
(\ref{-1.3}):
\begin{equation}
\label{2.1}
\|\psi(t)\| =: M(\psi) =: M 
\end{equation}
\begin{equation}
\label{2.2}
\|\nabla \psi(t)\|^2 + \frac{1}{2} (\|A^{\frac{1}{2}}\phi(t)\|^2 + \|\phi_t(t)\|^2) - \int_{{\bf R}^3} |\psi(t)|^2\phi(t)\, dx =: E =: E(\psi,\phi,\phi_t)
\end{equation}
By Gagliardo-Nirenberg we have
\begin{eqnarray}
\nonumber
\left| \int_{{\bf R}^3} |\psi|^2\phi dx\right| & \le & \|\phi\|_{L^6} \|\psi\|_{L^{\frac{12}{5}}}^2 \le c \|\nabla \phi\| \|\nabla \psi\|^{\frac{1}{2}} \|\psi\|^{\frac{3}{2}} \\
\label{2.3}
& \le & \frac{1}{4} \|A^{\frac{1}{2}}\phi\|^2 + \frac{1}{2} \|\nabla \psi\|^2 + c_1 \|\psi\|^6
\end{eqnarray}
This implies
\begin{equation}
\label{2.4}
\|\nabla \psi\|^2 + \frac{1}{2}(\|A^{\frac{1}{2}}\phi(t)\|^2 + \|\phi_t(t)\|^2) \le E + \frac{1}{4} \|A^{\frac{1}{2}}\phi(t)\|^2 + \frac{1}{2} \|\nabla \psi(t)\|^2 + c_1 M^6
\end{equation}
consequently
\begin{eqnarray}
\label{2.5}
\|\nabla \psi(t)\|^2 & \le & 2(E+c_1 M^6) \\
\label{2.6}
\|A^{\frac{1}{2}}\phi(t)\|^2 + \|\phi_t(t)\|^2 & \le & 4(E+c_1 M^6)
\end{eqnarray}
We also have
\begin{eqnarray}
\nonumber
E(\psi,\phi,\phi_t) & \le & \|\nabla \psi\|^2 + \frac{1}{2} (\|A^{\frac{1}{2}}\phi(t)\|^2 + \|\phi_t\|^2) + \left| \int |\psi|^2 \phi \, dx \right| \\
\label{2.7} & \le & \frac{3}{2} \|\nabla \psi\|^2 + \frac{3}{4} (\|A^{\frac{1}{2}}\phi(t)\|^2 + \|\phi_t\|^2) + c_1\|\psi\|^6
\end{eqnarray}
Consider now data
$$ \psi_0 \in H^{s,2}({\bf R}^3) \, , \, \phi_0 \in H^{m,2}({\bf R}^3) \, , \, \phi_1 \in H^{m-1,2}({\bf R}^3) $$
with $ 0 \le s,m \le 1 $.\\
We split these data into sums as follows:
$$ \psi_0 = \psi_{01} + \psi_{02} \, , \, \phi_0 = \phi_{01} + \phi_{02} \, , \, \phi_1 = \phi_{11} + \phi_{12} $$
where, for $N\ge 1$,
$$ \psi_{01}:= \int_{|\xi|\le N} e^{i<x,\xi>} \widehat{\psi_0}(\xi) \, d\xi $$
and $ \phi_{01} \, , \, \phi_{11} $ are defined in the same way.\\
One easily shows
\begin{eqnarray*}
\|\psi_{01}\|_{H^{l,2}} & \le & c N^{l-s} \|\psi_0\|_{H^{s,2}} \quad {\mbox for} \,\, l\ge s \quad , \quad \|\psi_{01} \| \le \|\psi_0\| \\
\|\psi_{02}\|_{H^{l,2}} & \le & c N^{l-s} \qquad {\mbox for} \,\, l \le s
\end{eqnarray*}
and similarly
\begin{eqnarray*}
\|\phi_{01}\|_{H^{l,2}} & \le & c N^{l-m} \|\phi_0\|_{H^{m,2}} \quad {\mbox for} \,\, l\ge m \quad , \quad \|\phi_{01} \| \le \|\phi_0\| \\
\|\phi_{02}\|_{H^{l,2}} & \le & c N^{l-m} \qquad {\mbox for} \,\, l \le m
\end{eqnarray*}
\begin{eqnarray*}
\|\phi_{11}\|_{H^{l-1,2}} & \le & c N^{l-m} \|\phi_1\|_{H^{m-1,2}} \quad {\mbox for} \,\, l\ge m  \\
\|\phi_{12}\|_{H^{l-1,2}} & \le & c N^{l-m} \qquad {\mbox for} \,\, l \le m
\end{eqnarray*}
Thus we have the following global bounds for the solution $(\widetilde{\psi},\widetilde{\phi})$ of ({\ref{-1.1}),(\ref{-1.2}) with data $(\psi_{01},\phi_{01},\phi_{11})$ - known to exist by \cite{B}, Thme. 3 - by (\ref{2.7}):
\begin{eqnarray}
\nonumber
E(\widetilde{\psi},\widetilde{\phi},\widetilde{\phi}_t) & \le & \frac{3}{2} \|\nabla \psi_{01}\|^2 + \frac{3}{4} (\|A^{\frac{1}{2}}\phi_{01}\|^2 + \|\phi_{11}\|^2) + c_1 \|\psi_{01}\|^6 \\
\label{2.7a}
& \le & \frac{{\overline c}}{2} ( N^{2(1-s)} + N^{2(1-m)}) \le \overline{c} N^{2(1-s\wedge m)}
\end{eqnarray}
and thus by (\ref{2.1}),(\ref{2.5}),(\ref{2.6}):
\begin{eqnarray}
\label{2.8}
\|\widetilde{\psi}(t)\| & \le &  M \\
\label{2.9}
\|\nabla \widetilde{\psi} (t) \| + \| A^{\frac{1}{2}} \widetilde{\phi}(t)\| + \|\widetilde{\phi}_t(t)\| & \le & \widehat{c} N^{1-s\wedge m}
\end{eqnarray}
The corresponding global solution $(\widetilde{\psi},\widetilde{\phi}_{\pm})$ of (\ref{1.4}),(\ref{1.5}) with data $ \psi_{01}$, $ \phi_{0\pm 1} := \phi_{01} \pm iA^{-\frac{1}{2}} \phi_{11} $ therefore fulfills
\begin{equation}
\label{2.10}
\|\nabla \widetilde{\psi}(t)\| + \|A^{\frac{1}{2}} \widetilde{\phi}_{\pm}(t)\| \le \widehat{c} N^{1-s\wedge m}
\end{equation}
where $\widehat{c}$ depends essentially only on $\overline{c}$ (the initial energy) and $M$ on the initial $L^2$-norm of $\widetilde{\psi}$.
\section{Further bounds for the regular part}
Consider the system of integral equations (\ref{1.6}),(\ref{1.7}) with $(\psi_0,\phi_{0\pm})$ replaced by $(\psi_{01},\phi_{0\pm 1})$ and $(\psi,\phi_{\pm})$ by $(\widetilde{\psi},\widetilde{\phi}_{\pm})$. Here $\phi_{0\pm1}:=\phi_{01} \pm iA^{-\frac{1}{2}} \phi_{11}$. \\
Let $ 0 \le \frac{2}{\gamma} = 3(\frac{1}{2} - \frac{1}{\rho}) < 1 $. By (\ref{-1.8}) we get
\begin{equation}
\label{3.2}
\|f\|_{X^{0,-\frac{1}{2}-}(I)} \le c \|f\|_{L^{\gamma '}_t(I,L^{\rho '}_x({\bf R}^3))}
\end{equation}
Interpolation with the trivial identity $ \|f\|_{X^{0,0}(I)} = \|f\|_{L^2_t(I,L^2_x({\bf R}^3))} $ gives
\begin{equation}
\label{3.3}
\|f\|_{X^{0,-\frac{1}{2}+}(I)} \le c \|f\|_{L^{\gamma '+}_t(I,L^{\rho '+}_x({\bf R}^3))}
\end{equation}
and also
\begin{equation}
\label{3.4}
\|f\|_{X^{1,-\frac{1}{2}+}(I)} \le c \|f\|_{L^{\gamma '+}_t(I,H^{1,\rho '+}_x({\bf R}^3))}
\end{equation}
We assume $|I| \le 1$. Applying this estimate to (\ref{1.6}) we get by (\ref{-1.5}) and the remarks following (\ref{1.7}):
\begin{eqnarray*}
\|\widetilde{\psi}\|_{X^{1,\frac{1}{2}+}(I)} & \le & c(\|\psi_{01}\|_{H^{1,2}} + \| \int^t_0 e^{i(t-s)\Delta} \frac{1}{2}(\widetilde{\phi}_+(s) + \widetilde{\phi}_-(s))\widetilde{\psi}(s)\, ds\|_{X^{1,\frac{1}{2}+}(I)}) \\
& \le & c(\|\psi_{01}\|_{H^{1,2}} + \|(\widetilde{\phi}_+ + \widetilde{\phi}_-)\widetilde{\psi}\|_{X^{1,-\frac{1}{2}+}(I)}) \\
 & \le & c(\|\psi_{01}\|_{H^{1,2}} + \|(\widetilde{\phi}_+ + \widetilde{\phi}_-)\widetilde{\psi}\|_{L^{\gamma '+}_t(I,H^{1,\rho '+}_x)})
\end{eqnarray*}
Choosing $ \gamma = 4-$ , $ \rho = 3+$ , $ \gamma '=\frac{4}{3}+$ , $ \rho ' = \frac{3}{2}-$ we estimate
\begin{eqnarray*}
\|\widetilde{\phi}_{\pm} \widetilde{\psi}\|_{H^{1,\rho '+}_x} & \le & c(\|(\nabla \widetilde{\phi}_{\pm})\widetilde{\psi}\|_{L^{\rho '+}_x} + \| \widetilde{\phi}_{\pm} \nabla \widetilde{\psi}\|_{L^{\rho '+}_x}) + \| \widetilde{\phi}_{\pm} \widetilde{\psi}\|_{L^{\rho '+}_x}) \\
& \le & c(\|\nabla \widetilde{\phi}_{\pm}\|_{L^2_x}\|\widetilde{\psi}\|_{L^6_x}  + \| \widetilde{\phi}_{\pm}\|_{L^6_x} \| \nabla \widetilde{\psi}\|_{L^2_x} + \| \widetilde{\phi}_{\pm}\|_{L^6_x} \| \widetilde{\psi}\|_{L^2_x}) \\
& \le & c \|\nabla \widetilde{\phi}_{\pm}\|_{L^2_x} (\|\nabla \widetilde{\psi}\|_{L^2_x} + \|\widetilde{\psi}\|_{L^2_x}) \\
& \le & c N^{2(1-s\wedge m)}
\end{eqnarray*}
Thus
$$ \|\widetilde{\psi}\|_{X^{1,\frac{1}{2}+}(I)} \le c (N^{1-s} + N^{2(1-s\wedge m)} |I|^{\frac{3}{4}-})$$
Assume now $ |I| \le N^{-\frac{4}{3}(1-s\wedge m)-} $ . Then we conclude
$$\|\widetilde{\psi}\|_{X^{1,\frac{1}{2}+}(I)} \le c N^{1-s\wedge m} $$
The same estimate also holds true, if we only assume $\|\psi_{01}\|_{H^{1,2}} \le cN^{1-s\wedge m} $ (important remark for the iteration process described later).\\
Next we estimate $\|\widetilde{\psi}\|_{X^{0,\frac{1}{2}+}(I)}$. We again use (\ref{3.3}) with $\gamma =4-$, $\rho = 3+$ and conclude as before
\begin{eqnarray*}
\|\widetilde{\phi}_{\pm} \widetilde{\psi}\|_{L^{\gamma '+}_t(I,L^{\rho '+}_x)} & \le & \|\widetilde{\psi}\|_{L^{\infty}_t(I,L^2_x)} \|\widetilde{\phi}_{\pm}\|_{L^{\infty}_t(I,L^6_x)} |I|^{\frac{3}{4}-} \\
& \le & cMN^{1-s\wedge m} N^{-(1-s\wedge m)} \\
& \le & c
\end{eqnarray*}
Applying this estimate to ({\ref{1.6}) we get as before
$$ \|\widetilde{\psi}\|_{X^{0,\frac{1}{2}+}(I)} \le c (\|\psi_{01}\|_{L^2_x} +1 ) \le c $$
if $ \|\psi_{01}\|_{L^2} \le c $.\\
In order to estimate $\|\widetilde{\phi}_{\pm}\|_{X^{1,\frac{1}{2}+}(I)}$ we start with (\ref{-1.10}) and get $\|f\|_{X^{0,-\frac{1}{2}-}_{\pm}(I)}$\\ $ \le c \|f\|_{L^1_t(I,L^2_x)} $. Interpolation with the trivial identity $ \|f\|_{X^{0,0}_{\pm}(I)} = \|f\|_{L^2_t(I,L^2_x)}$ gives
\begin{equation}
\label{3.5}
\|f\|_{X^{0,-\frac{1}{2}+}_{\pm}(I)} \le c \|f\|_{L^{1+}_t(I,L^2_x)}
\end{equation}
Applying this estimate to (\ref{1.7}) gives
\begin{eqnarray*}
\|\widetilde{\phi}_{\pm}\|_{X^{1,\frac{1}{2}+}_{\pm}(I)} & \le & c(\|\phi_{0\pm 1}\|_{H^{1,2}_x} + \| \int_0^t e^{\mp i(t-s)A^{\frac{1}{2}}} A^{-\frac{1}{2}}(|\widetilde{\psi}(s)|^2)ds\|_{X^{1,\frac{1}{2}+}_{\pm}(I)}) \\
& \le & c(\|\phi_{0\pm 1}\|_{H^{1,2}_x} + \| |\widetilde{\psi}|^2 \|_{X^{0,-\frac{1}{2}+}_{\pm}(I)}) \\
& \le & c(\|\phi_{0\pm 1}\|_{H^{1,2}_x} + \| |\widetilde{\psi}|^2 \|_{L^{1+}_t(I,L^2_x)})
\end{eqnarray*}
We have
$$ \| |\widetilde{\psi}|^2\|_{L^2_x} \le c \|\widetilde{\psi}\|_{L^2_x}^{\frac{1}{2}} \|\widetilde{\psi}\|_{L^6_x}^{\frac{3}{2}} \le c \|\widetilde{\psi}\|_{L^2_x}^{\frac{1}{2}} \|\nabla \widetilde{\psi}\|_{L^2_x}^{\frac{3}{2}} $$
and consequently
\begin{eqnarray*}
\|\widetilde{\phi}_{\pm}\|_{X^{1,\frac{1}{2}+}_{\pm}(I)} & \le & c (\|\phi_{0\pm 1}\|_{H^{1,2}_x} + \|\widetilde{\psi}\|_{L^{\infty}_t(I,L^2_x)}^{\frac{1}{2}} \|\nabla \widetilde{\psi}\|_{L^{\infty}_t(I,L^2_x)}^{\frac{3}{2}} |I|^{1-} ) \\
& \le & c(\|\phi_{0\pm 1}\|_{H^{1,2}_x} + M^{\frac{1}{2}} N^{\frac{3}{2}(1-s\wedge m)} N^{-\frac{4}{3}(1-s\wedge m)}) \\
& \le & c (N^{1-m} + N^{\frac{1}{6}(1-s\wedge m)})
\end{eqnarray*}
Especially we conclude
$$ \|\widetilde{\phi}_{\pm}\|_{X^{1,\frac{1}{2}+}_{\pm}(I)} \le c N^{1-s\wedge m} $$
if $ \|\phi_{0\pm 1}\|_{H^{1,2}_x} \le c N^{1-s\wedge m} $ , i.e. $ \|\phi_{01}\|_{H_x^{1,2}} + \|\phi_{11}\|_{L_x^2} \le c N^{1-s\wedge m}$.\\[2ex]
We summarize the results obtained so far in the following 
\begin{lemma}
\label{Lemma 3.1}
If $|I|\le N^{-\frac{4}{3}(1-s\wedge m)-}$ and
\begin{equation}
\label{3.5a}
\|\psi_{01}\|_{L^2_x} \le c \quad , \quad \|\psi_{01}\|_{H^{1,2}_x} + \|\phi_{01}\|_{H^{1,2}_x} + \|\phi_{11}\|_{L^2_x} \le c N^{1-s\wedge m}
\end{equation}
i.e.
\begin{equation}
\label{3.6}
\|\psi_{01}\|_{L^2_x} \le c \quad , \quad \|\psi_{01}\|_{H^{1,2}_x} + \|\phi_{0\pm 1}\|_{L^2_x}  \le c N^{1-s\wedge m}
\end{equation}
the following estimates hold:
$$
\|\widetilde{\psi}\|_{X^{0,\frac{1}{2}+}(I)} \le c \quad , \quad 
\|\widetilde{\psi}\|_{X^{1,\frac{1}{2}+}(I)}  \le  c N^{1-s\wedge m} \quad ,  \quad 
\|\widetilde{\phi}_{\pm}\|_{X_{\pm}^{1,\frac{1}{2}+}(I)} \le c N^{1-s\wedge m}
$$
Also the estimates (\ref{2.8}),(\ref{2.9}),(\ref{2.10}) hold under these assumptions.
\end{lemma}
{\bf Important remark:} Here and in the sequel the constants denoted by $c$ depend essentially only on $\overline{c}$ in (\ref{2.7a}) (and therefore on $E(\widetilde{\psi},\widetilde{\phi},\widetilde{\phi}_t)$ and on $M$).

\section{The part with rough data}
Let $(\psi,\phi_+,\phi_-)$ be a solution of (\ref{1.4}),(\ref{1.5}) with data $(\psi_0,\phi_{0+},\phi_{0-})$ and  $(\widetilde{\psi},\widetilde{\phi}_+,\widetilde{\phi}_-)$ be the solution with data $(\psi_{01},\phi_{0+1},\phi_{0-1})$.\\
Define $ \widehat{\psi} := \psi - \widetilde{\psi} $ , $ \widehat{\phi}_{\pm} := \phi_{\pm} - \widetilde{\phi}_{\pm} $. Then $(\widehat{\psi},\widehat{\phi}_+,\widehat{\phi}_-)$ fulfills
\begin{eqnarray}
\nonumber
i\widehat{\psi}_t + \Delta \widehat{\psi} & = & i \psi_t + \Delta \psi -i \widetilde{\psi}_t + \Delta \widetilde{\psi} = -\psi \phi + \widetilde{\psi} \widetilde{\phi} \\
\nonumber
& = & - \frac{1}{2} \psi (\phi_+ + \phi_-) + \frac{1}{2}\widetilde{\psi} (\widetilde{\phi}_+ + \widetilde{\phi}_-) \\
\label{4.0a}
& = & - \frac{1}{2} (\widehat{\psi} + \widetilde{\psi}) (\widehat{\phi}_+ + \widehat{\phi}_- + \widetilde{\phi}_+ + \widetilde{\phi}_-) + \frac{1}{2}\widetilde{\psi} (\widetilde{\phi}_+ + \widetilde{\phi}_-) \\
\nonumber
& = & -\frac{1}{2} \widehat{\psi} (\widetilde{\phi}_+ + \widetilde{\phi}_-) - \frac{1}{2} \widehat{\psi}(\widehat{\phi}_+ + \widehat{\phi}_-) - \frac{1}{2} \widetilde{\psi}(\widehat{\phi}_+ + \widehat{\phi}_-) \\
\nonumber
& =: & F_1 + F_2 + F_3 =: F
\end{eqnarray}
and
\begin{eqnarray}
\nonumber
i\widehat{\phi}_{\pm t} \mp A^{1/2} \widehat{\phi}_{\pm} & = & i\phi_{\pm t} \mp A^{1/2} \phi_{\pm} - i\widetilde{\phi}_{\pm t} \pm A^{1/2} \widetilde{\phi}_{\pm} \\
\nonumber
& = & \pm A^{-1/2} (|\psi|^2) \mp A^{-1/2} (|\widetilde{\psi}|^2) \\
\label{4.0b}
& = & \pm A^{-1/2} ((\widehat{\psi} + \widetilde{\psi})(\overline{\widehat{\psi}} + \overline{\widetilde{\psi}}) \mp A^{-1/2} (\widetilde{\psi} \overline{\widetilde{\psi}})  \\
\nonumber
& = & \pm A^{-1/2} (|\widehat{\psi}|^2) \pm A^{-1/2} (\overline{\widehat{\psi}} \widetilde{\psi}) \pm A^{-1/2} (\widehat{\psi} \overline{\widetilde{\psi}}) \\
\nonumber
& =: & G_1 + G_2 + G_3 =: G
\end{eqnarray} 
Furthermore

\begin{eqnarray*}
\widehat{\psi}(0) & = & \psi(0) - \widetilde{\psi}(0) = \psi_0 - \psi_{01} = \psi_{02} \\
\widehat{\phi}_{\pm}(0) & = & \phi_{\pm}(0) - \widetilde{\phi}_{\pm}(0) = \phi_{0\pm} - \phi_{0\pm 1} =: \phi_{0\pm 2}
\end{eqnarray*}
The corresponding system of integral equations reads as follows:
\begin{eqnarray}
\label{4.1}
\widehat{\psi}(t) & = & e^{it\Delta} \psi_{02} - i \int^t_0 e^{i(t-s)\Delta} F(s) ds =: e^{it\Delta} \psi_{02} + w(t) \\
\label{4.2}
\widehat{\phi}_{\pm}(t) & = & e^{\mp itA^{1/2}} \phi_{0\pm 2} - i \int^t_0 e^{\mp i(t-s)A^{1/2}} G(s) ds =: e^{\mp itA^{1/2}} \phi_{0\pm 2} + z_{\pm}(t)
\end{eqnarray}  
Here we have
\begin{eqnarray}
\label{4.3}
\|\psi_{02}\|_{H^{s,2}_x} \le c & , & \|\psi_{02}\|_{L^2_x} \le c N^{-s} \\
\label{4.4}
\|\phi_{0\pm 2}\|_{H^{m,2}_x} \le c & , & \|\phi_{0\pm 2}\|_{L^2_x} \le c N^{-m}
\end{eqnarray}
We construct a solution of (\ref{4.1}),(\ref{4.2}) in some time interval $I$ using the contraction mapping principle.\\
We define a mapping $S = (S_0,S_+,S_-)$ by
\begin{eqnarray*}
(S_0 \widehat{\psi})(t) & := & e^{it\Delta} \psi_{02} + w(t) \\
(S_{\pm} \widehat{\phi}_{\pm})(t) & := & e^{\mp itA^{1/2}} \phi_{0\pm 2} + z_{\pm}(t)
\end{eqnarray*}
\begin{prop}
\label{Proposition 4.1}
For $ 1\ge s,m > \frac{7}{10} $ and data $\psi_{02} \in H^{s,2}_x({\bf R}^3)$, $ \phi_{0\pm 2} \in H^{m,2}_x({\bf R}^3)$ with (\ref{4.3}),(\ref{4.4}) and $\psi_{01}, \phi_{0\pm 1}$ as in (\ref{3.5a}),(\ref{3.6}) the system of integral equations (\ref{4.1}),(\ref{4.2}) has a unique solution $(\widehat{\psi},\widehat{\phi}_{\pm}) \in X^{s,\frac{1}{2}+}(I) \times X^{m,\frac{1}{2}+}_{\pm}(I)$ in the same interval $I$ of the preceding section with $|I| \le N^{-\frac{4}{3}(1-s\wedge m)-\delta}$ $(\delta > 0)$, which fulfills
\begin{enumerate}
\item in the case $ s\le m$
\begin{eqnarray*}
\|\widehat{\psi}\|_{X^{0,\frac{1}{2}+}(I)} & \le & c N^{-s} \le c N^{-s\wedge m} \\
\|\widehat{\phi}_{\pm}\|_{X^{0,\frac{1}{2}+}_{\pm}(I)} & \le & c (N^{-\frac{5}{6}-\frac{1}{6}s} + N^{-m} ) \le c N^{-s\wedge m} \\
\end{eqnarray*}
\item in the case $ m\le s$
\begin{eqnarray*}
\|\widehat{\psi}\|_{X^{0,\frac{1}{2}+}(I)} & \le & c (N^{-\frac{1}{6}-\frac{5}{6}m}+N^{-s})  \le c N^{-s\wedge m} \\
\|\widehat{\phi}_{\pm}\|_{X^{0,\frac{1}{2}+}_{\pm}(I)} & \le & c N^{-m}  \le c N^{-s\wedge m} \\
\end{eqnarray*}
\end{enumerate}
and moreover in any case 
$$ \|\widehat{\psi}\|_{X^{s,\frac{1}{2}+}(I)} \le c \quad , \quad
\|\widehat{\phi}\|_{X_{\pm}^{m,\frac{1}{2}+}(I)} \le c $$
\end{prop}
{\bf Proof:}
We use Banach's fixed point theorem in the set $Z$, where
\begin{enumerate}
\item in the case $s\le m$
\begin{eqnarray*}
Z & := & \{ \|\widehat{\psi}\|_{X^{s,\frac{1}{2}+}(I)} \le c_0 \, , \, \|\widehat{\psi}\|_{X^{0,\frac{1}{2}+}(I)} \le c_0 N^{-s} \\
& & \, \, \|\widehat{\phi}_{\pm}\|_{X^{m,\frac{1}{2}+}_{\pm}(I)} \le c_0 \, , \, \|\widehat{\phi}_{\pm}\|_{X^{0,\frac{1}{2}+}_{\pm}(I)} \le c_0 (N^{-\frac{5}{6}-\frac{1}{6}s} + N^{-m}) \}
\end{eqnarray*}
\item in the case $m\le s$
\begin{eqnarray*}
Z & := & \{ \|\widehat{\psi}\|_{X^{s,\frac{1}{2}+}(I)} \le c_0 \, , \, \|\widehat{\psi}\|_{X^{0,\frac{1}{2}+}(I)} \le c_0 (N^{-\frac{1}{6}-\frac{5}{6}m} + N^{-s})  \\
& & \, \, \|\widehat{\phi}_{\pm}\|_{X^{m,\frac{1}{2}+}_{\pm}(I)} \le c_0 \, , \, \|\widehat{\phi}_{\pm}\|_{X^{0,\frac{1}{2}+}_{\pm}(I)} \le c_0 N^{-m} \}
\end{eqnarray*}
$c_0$ is chosen below.
\end{enumerate}
Now take any $(\widehat{\psi},\widehat{\phi}_+,\widehat{\phi}_-) \in Z$. In order to show $(S_0\widehat{\psi},S_+\widehat{\phi}_+,S_-\widehat{\phi}_-) \in Z$ we first estimate $ \|S_0\widehat{\psi}\|_{X^{s,\frac{1}{2}+}(I)}$.\\
By Lemma \ref{Lemma 0.5} with $\sigma = 0+$ we get by interpolation
\begin{eqnarray}
\nonumber
\|\widehat{\psi}\widehat{\phi}_{\pm}\|_{X^{s+,-\frac{1}{2}-}(I)}  
&  \le & c(\|\widehat{\psi}\|_{X^{0,\frac{1}{2}+}(I)} \|\widehat{\phi}_{\pm}\|_{X^{s-\frac{1}{2}+,\frac{1}{2}+}_{\pm}(I)} +
\|\widehat{\psi}\|_{X^{s,0}(I)} \|\widehat{\phi}_{\pm}\|_{X^{\frac{1}{2}+,\frac{1}{2}+}_{\pm}(I)}) \\
\label{4.5}
&  \le & c(\|\widehat{\psi}\|_{X^{0,\frac{1}{2}+}(I)} \|\widehat{\phi}_{\pm}\|_{X^{m,\frac{1}{2}+}_{\pm}(I)}^{\frac{s-\frac{1}{2}}{m}+} \|\widehat{\phi}_{\pm}\|_{X^{0,\frac{1}{2}+}_{\pm}(I)}^{1-\frac{s-\frac{1}{2}}{m}-} \\ 
\nonumber
& & \quad +
\|\widehat{\psi}\|_{X^{s,\frac{1}{2}+}(I)} |I|^{\frac{1}{2}} \|\widehat{\phi}_{\pm}\|_{X^{m,\frac{1}{2}+}_{\pm}(I)}^{\frac{1}{2m}+}
\|\widehat{\phi}_{\pm}\|_{X^{0,\frac{1}{2}+}_{\pm}(I)}^{1-\frac{1}{2m}-}) \\
\nonumber
& \le & c(N^{-s\wedge m} N^{-s\wedge m(1-\frac{s-\frac{1}{2}}{m})+} + N^{-\frac{2}{3}(1-s\wedge m)-\frac{\delta}{2}} N^{-s\wedge m(1-\frac{1}{2m}-)}) \\
\nonumber
& \le & c N^{-\gamma(s,m)}
\end{eqnarray}
with $ \gamma(s,m) > 0 $ for $ 1\ge s,m>\frac{1}{2}$.\\
Next, by Lemma \ref{Lemma 0.5} with $\sigma = \frac{1}{2}-$ and Lemma \ref{Lemma 3.1}:
\begin{eqnarray}
\nonumber
\lefteqn{\|\widehat{\psi}\widetilde{\phi}_{\pm}\|_{X^{s+,-\frac{1}{2}-}(I)}} \\ \nonumber
& & \le c(\|\widehat{\psi}\|_{X^{0,\frac{1}{2}+}(I)} \|\widetilde{\phi}_{\pm}\|_{X^{s-\frac{1}{2}+,\frac{1}{2}+}_{\pm}(I)} +
\|\widehat{\psi}\|_{X^{s-\frac{1}{2}+,0}(I)} \|\widetilde{\phi}_{\pm}\|_{X^{1,\frac{1}{2}+}_{\pm}(I)}) \\
\nonumber
& & \le  c(\|\widehat{\psi}\|_{X^{0,\frac{1}{2}+}(I)} \|\widetilde{\phi}_{\pm}\|_{X^{s-\frac{1}{2}+,\frac{1}{2}+}_{\pm}(I)} +
\|\widehat{\psi}\|_{X^{s-\frac{1}{2}+,\frac{1}{2}+}(I)} |I|^{\frac{1}{2}} \|\widetilde{\phi}_{\pm}\|_{X^{1,\frac{1}{2}+}_{\pm}(I)}) \\
\nonumber
& & \le  c(\|\widehat{\psi}\|_{X^{0,\frac{1}{2}+}(I)} \|\widetilde{\phi}_{\pm}\|_{X^{1,\frac{1}{2}+}_{\pm}(I)} +
\|\widehat{\psi}\|_{X^{s,\frac{1}{2}+}(I)}^{\frac{s-\frac{1}{2}}{s}+} \|\widehat{\psi}\|_{X^{0,\frac{1}{2}+}(I)}^{1-\frac{s-\frac{1}{2}}{s}-} |I|^{\frac{1}{2}} \|\widetilde{\phi}_{\pm}\|_{X^{1,\frac{1}{2}+}_{\pm}(I)}) \\
\nonumber
& & \le  c(N^{-s\wedge m} N^{1-s\wedge m} + N^{-s\wedge m(1-\frac{s-\frac{1}{2}}{s})-} N^{-\frac{2}{3}(1-s\wedge m)-\frac{\delta}{2}} N^{1-s\wedge m}) \\
\label{4.6}
& & \le  c (N^{1-2(s\wedge m)} + N^{\frac{1}{3}-s\wedge m(\frac{1}{2s}+\frac{1}{3})}) \\
\nonumber
& & \le  c (N^{1-2(s\wedge m)} + N^{\frac{1}{3}-\frac{1}{2}(\frac{1}{2}+\frac{1}{3})}) \\
\nonumber
& & \le  c N^{-\gamma(s,m)}
\end{eqnarray}	
Similarly by Lemma \ref{Lemma 0.5} with $\sigma = 0+$ and Lemma \ref{Lemma 3.1}:
\begin{eqnarray*}
\lefteqn{\|\widetilde{\psi}\widehat{\phi}_{\pm}\|_{X^{s+,-\frac{1}{2}-}(I)}} \\ 
& & \le c(\|\widetilde{\psi}\|_{X^{0,\frac{1}{2}+}(I)} \|\widehat{\phi}_{\pm}\|_{X^{s-\frac{1}{2}+,\frac{1}{2}+}_{\pm}(I)} +
\|\widetilde{\psi}\|_{X^{s,0}(I)} \|\widehat{\phi}_{\pm}\|_{X^{\frac{1}{2}+,\frac{1}{2}+}_{\pm}(I)}) \\
& & \le c(\|\widetilde{\psi}\|_{X^{0,\frac{1}{2}+}(I)} \|\widehat{\phi}_{\pm}\|_{X^{m,\frac{1}{2}+}_{\pm}(I)}^{\frac{s-\frac{1}{2}}{m}+} \|\widehat{\phi}_{\pm}\|_{X^{0,\frac{1}{2}+}_{\pm}(I)}^{1-\frac{s-\frac{1}{2}}{m}-} \\
& & \quad +
\|\widetilde{\psi}\|_{X^{0,\frac{1}{2}+}(I)}^{1-s}
\|\widetilde{\psi}\|_{X^{1,\frac{1}{2}+}(I)}^s |I|^{\frac{1}{2}} \|\widehat{\phi}_{\pm}\|_{X^{m,\frac{1}{2}+}_{\pm}(I)}^{\frac{1}{2m}+} \|\widehat{\phi}_{\pm}\|_{X^{0,\frac{1}{2}+}_{\pm}(I)}^{1-\frac{1}{2m}-}) \\
& & \le c(N^{-(s\wedge m)(1-\frac{s-\frac{1}{2}}{m})+} + N^{s(1-s\wedge m)} N^{-\frac{2}{3}(1-s\wedge m)-\frac{\delta}{2}} N^{-(s\wedge m)(1-\frac{1}{2m}-)} ) \\
& & \le c (N^{-(s\wedge m)(1-\frac{s-\frac{1}{2}}{m})+} + N^{s-\frac{2}{3}-(s\wedge m)(s+\frac{1}{3}-\frac{1}{2m})+})
\end{eqnarray*}
Now, for $ s,m > \frac{7}{10} $ the last exponent is \\
$ \le s-\frac{2}{3}-\frac{7}{10}(s+\frac{1}{3}-\frac{5}{7}) = \frac{3}{10}s - \frac{2}{3} - \frac{7}{30} + \frac{1}{2} = \frac{3}{10}s - \frac{2}{5} \le \frac{3}{10}-\frac{2}{5} = - \frac{1}{10} < 0 $ \\
Thus
\begin{equation}
\label{4.7}
\|\widetilde{\psi} \widehat{\phi}_{\pm} \|_{X^{s+,-\frac{1}{2}-}(I)} \le c N^{-\gamma(s,m)}
\end{equation}
Now we interpolate between (\ref{4.6}) and the following estimate, which follows from (\ref{1.10}) and Lemma \ref{Lemma 3.1}:
\begin{eqnarray*}
\|\widehat{\psi}\widetilde{\phi}_{\pm}\|_{X^{0,-\frac{1}{4}-\epsilon}(I)} & \le & c(\|\widetilde{\phi}_{\pm}\|_{X^{0,\frac{1}{2}+\epsilon}_{\pm}(I)} \|\widehat{\psi}\|_{X^{0,\frac{1}{2}+\epsilon}_{\pm}(I)}) \\
& \le & c N^{1-s\wedge m} N^{-s\wedge m} = c N^{1-2(s\wedge m)}
\end{eqnarray*}
This gives
\begin{equation}
\label{4.8}
\|\widehat{\psi}\widetilde{\phi}_{\pm}\|_{X^{s,-\frac{1}{2}+}(I)} \le c N^{-\gamma(s,m)+}
\end{equation}
Similarly we treat (\ref{4.5}) and (\ref{4.7}) and conclude
$$ \|F\|_{X^{s,-\frac{1}{2}+}(I)} \le c N^{-\gamma(s,m)+} $$
Therefore, with $ c_0 \ge 2c\|\psi_{02}\|_{H^{s,2}({\bf R}^3)} $ , we have
$$
\|S_0 \widehat{\psi}\|_{X^{s,\frac{1}{2}+}(I)} \le c \|\psi_{02}\|_{H^{s,2}} + c \|F\|_{X^{s,-\frac{1}{2}+}(I)} 
 \le  \frac{c_0}{2} + c N^{-\gamma(s,m)+} \le c_0
$$
if $N$ is sufficiently large.\\
Next we estimate $ \|S_0 \widehat{\psi}\|_{X^{0,\frac{1}{2}+}(I)} $. Using (\ref{4.1}) we have to control $\|F\|_{X^{0,-\frac{1}{2}+}(I)}$. By (\ref{3.3}) with $\gamma ' = \frac{4}{3}+ $ , $ \rho ' = \frac{3}{2}- $ and Lemma \ref{Lemma 3.1} we have
\begin{eqnarray}
\nonumber
\|\widehat{\psi} \widetilde{\phi}_{\pm}\|_{X^{0,-\frac{1}{2}+}(I)} & \le & c \|\widehat{\psi} \widetilde{\phi}_{\pm}\|_{L^{\frac{4}{3}+}_t(I,L^{\frac{3}{2}}_x)} \\ \nonumber
& \le & c \|\widehat{\psi}\|_{L^{\infty}_t(I,L^2_x)} \|\widetilde{\phi}_{\pm}\|_{L^{\infty}_t(I,L^6_x)} |I|^{\frac{3}{4}-} \\
\label{4.9}
& \le & c N^{1-(s\wedge m)} N^{-(1-s\wedge m)-\frac{3}{4}\delta +} \|\widehat{\psi}\|_{X^{0,\frac{1}{2}+}(I)} \\ \nonumber
& \le & c N^{-\frac{3}{4}\delta +} \|\widehat{\psi}\|_{X^{0,\frac{1}{2}+}(I)}
\end{eqnarray}
Similarly by (\ref{3.3}) with $\gamma ' =2$ , $ \rho ' = \frac{6}{5} $ we get
\begin{eqnarray}
\nonumber
\|\widehat{\psi} \widehat{\phi}_{\pm}\|_{X^{0,-\frac{1}{2}+}(I)} & \le & c \|\widehat{\psi} \widehat{\phi}_{\pm}\|_{L^{2+}_t(I,L^{\frac{6}{5}+}_x)} \\ \nonumber
& \le & c \|\widehat{\psi}\|_{L^{\infty}_t(I,L^2_x)} \|\widehat{\phi}_{\pm}\|_{L^{\infty}_t(I,L^{3+}_x)} |I|^{\frac{1}{2}-} \\ \nonumber
& \le & c \|\widehat{\psi}\|_{L^{\infty}_t(I,L^2_x)} \|\widehat{\phi}_{\pm}\|_{L^{\infty}_t(I,L^2_x)}^{1-\frac{1}{2m}-} \|\widehat{\phi}_{\pm}\|_{L^{\infty}_t(I,H^{m,2}_x)}^{\frac{1}{2m}+} |I|^{\frac{1}{2}-} \\
\label{4.10}
& \le & c N^{-(s\wedge m)(1-\frac{1}{2m})+} N^{-\frac{2}{3}(1-s\wedge m)-\frac{
\delta}{2} +} \|\widehat{\psi}\|_{X^{0,\frac{1}{2}+}(I)} \\ \nonumber
& \le & c N^{-\frac{\delta}{2} +} \|\widehat{\psi}\|_{X^{0,\frac{1}{2}+}(I)}
\end{eqnarray}
Finally by (\ref{3.3}) with $ \gamma ' =2$ , $ \rho ' = \frac{6}{5} $ again and Lemma \ref{Lemma 3.1} we get
\begin{eqnarray}
\nonumber
\|\widetilde{\psi} \widehat{\phi}_{\pm}\|_{X^{0,-\frac{1}{2}+}(I)} & \le & c \|\widetilde{\psi} \widehat{\phi}_{\pm}\|_{L^{2+}_t(I,L^{\frac{6}{5}+}_x)} \\ \nonumber
& \le & c \|\widetilde{\psi}\|_{L^{\infty}_t(I,L^{3+}_x)} \|\widehat{\phi}_{\pm}\|_{L^{\infty}_t(I,L^2_x)} |I|^{\frac{1}{2}-} \\ \nonumber
& \le & c \|\widetilde{\psi}\|_{L^{\infty}_t(I,L^2_x)}^{\frac{1}{2}-} \|\nabla \widetilde{\psi}\|_{L^{\infty}_t(I,L^2_x)}^{\frac{1}{2}+}\|\widehat{\phi}_{\pm}\|_{L^{\infty}_t(I,L^2_x)} |I|^{\frac{1}{2}-} \\
\label{4.11}
& \le & c N^{\frac{1}{2}(1-s\wedge m)+} N^{-\frac{2}{3}(1-s\wedge m)-\frac{
\delta}{2} +} \|\widehat{\phi}_{\pm}\|_{X_{\pm}^{0,\frac{1}{2}+}(I)} \\ \nonumber
& = & c N^{-\frac{1}{6}(1-s\wedge m)-\frac{\delta}{2} +} \|\widehat{\phi}_{\pm}\|_{X^{0,\frac{1}{2}+}_{\pm}(I)}
\end{eqnarray}
\begin{enumerate}
\item If now $s\le m$ we conclude from (\ref{4.11})
\begin{eqnarray*}
\|\widetilde{\psi}\widehat{\phi}_{\pm}\|_{X^{0,-\frac{1}{2}+}(I)} & \le & c N^{-\frac{1}{6}(1-s)-\frac{\delta}{2}+} (N^{-\frac{5}{6}-\frac{1}{6}s} + N^{-m}) \\
& = & c N^{-1-\frac{\delta}{2}+} + N^{-m-\frac{1}{6}(1-s)-\frac{\delta}{2}+}
\end{eqnarray*}
Now $ -m-\frac{1}{6}(1-s) \le -s \Leftrightarrow \frac{7}{6}s \le m + \frac{1}{6} $. This is fulfilled because $ \frac{7}{6}s \le \frac{7}{6}m = m + \frac{1}{6}m \le m + \frac{1}{6} $. Thus
\begin{equation}
\label{4.12}
\|\widetilde{\psi}\widehat{\phi}_{\pm}\|_{X^{0,-\frac{1}{2}+}(I)} \le c N^{-s-}
\end{equation}
\item If $m\le s$, we conclude from (\ref{4.11})
\begin{equation}
\label{4.13}
\|\widetilde{\psi}\widehat{\phi}_{\pm}\|_{X^{0,-\frac{1}{2}+}(I)} \le c N^{-\frac{1}{6}(1-m)-\frac{\delta}{2}+} N^{-m} = c N^{-\frac{1}{6} - \frac{5}{6}m - \frac{\delta}{2}+}
\end{equation}
\end{enumerate}
From (\ref{4.1}),(\ref{4.9}),(\ref{4.10}),(\ref{4.12}),(\ref{4.14}) we conclude
\begin{enumerate}
\item in the case $s\le m$
\begin{eqnarray}
\nonumber
\|S_0 \widehat{\psi}\|_{X^{0,\frac{1}{2}+}(I)} & \le & c \|\psi_{02}\|_{L^2} + c N^{-\frac{\delta}{2}+} \|\widehat{\psi}\|_{X^{0,\frac{1}{2}+}(I)} + c N^{-s-} \\
\nonumber & \le & c N^{-s} + c N^{-s-} \\
\label{4.14}
& \le & \frac{c_0}{2} N^{-s} + c N^{-s-} \\ \nonumber
& \le & c_0 N^{-s}
\end{eqnarray}
if $c_0 \ge 2c$ and $N$ sufficiently large.
\item Similarly for $m\le s$ we get
\begin{eqnarray}
\nonumber
\|S_0 \widehat{\psi}\|_{X^{0,\frac{1}{2}+}(I)} & \le & c N^{-s} + c N^{-s-} + c N^{-\frac{1}{6} - \frac{5}{6}m - \frac{\delta}{2}+} \\
\label{4.15}
& \le & c_0 (N^{-s} + N^{-\frac{1}{6} - \frac{5}{6}m})
\end{eqnarray}
\end{enumerate}
Next we estimate $ \|S_{\pm}\widehat{\phi}_{\pm}\|_{X^{m,\frac{1}{2}+}_{\pm}(I)}$. Assume
\begin{equation}
\label{4.16}
0 \le \frac{2}{\gamma} = 2(\frac{1}{2}-\frac{1}{\rho}) \quad {\mbox and} \quad l+2(\frac{1}{2}-\frac{1}{\rho}) = 0
\end{equation}
Then we get by (\ref{-1.11}):
$$ \|f\|_{X^{0,-\frac{1}{2}-}_{\pm}} \le c \|f\|_{L^{\gamma '}_t({\bf R},H^{-l,\rho '}_x)} $$
Interpolating with the trivial identity $\|f\|_{X^{0,0}_{\pm}} = \|f\|_{L^2_t({\bf R},L^2_x)}$ we have
$$ \|f\|_{X^{0,-\frac{1}{2}+}_{\pm}} \le c \|f\|_{L^{\gamma '+}_t({\bf R},H^{-l+,\rho '+}_x)} $$
and also
\begin{equation}
\label{4.17}
\|f\|_{X^{m-1,-\frac{1}{2}+}_{\pm}} \le c \|f\|_{L^{\gamma '+}_t({\bf R},H^{m-1-l+,\rho '+}_x)}
\end{equation}
In order to estimate $G_2$ we use (\ref{4.16}) with $l=0$ , $ \rho =2$ , $ \gamma = \infty $ and get
$$ \|\overline{\widehat{\psi}}\widetilde{\psi}\|_{L^{1+}(I,H^{m-1,2})} \le c \|\overline{\widehat{\psi}}\widetilde{\psi}\|_{L^{1+}(I,L^{\frac{6}{5-2m}})} \le c\left(\int_I \|\widehat{\psi}\|^{1+}_{L^{\frac{6}{5-2m}\widehat{p}}} \|\widetilde{\psi}\|^{1+}_{L^{\frac{6}{5-2m}\widehat{q}}}\, dt\right)^{1-} $$
Choosing $\frac{1}{\widehat{p}} = \frac{3-2s}{5-2m} $ ($\le 1$, because $m-s\le 1$) and $\frac{1}{\widehat{q}} = \frac{2(1+s-m)}{5-2m} $ we have $H^{s,2} \subset L^{\frac{6}{5-2m}\widehat{p}}$ and $ \frac{6}{5-2m}\widehat{q} = \frac{3}{1+s-m}$, thus by interpolation
$$ \|\widetilde{\psi}\|_{L^{\frac{6\widehat{q}}{5-2m}}} \le c \|\widetilde{\psi}\|_{L^2}^{s-m+\frac{1}{2}} \|\widetilde{\psi}\|_{H^{1,2}}^{\frac{1}{2}+m-s} $$
We get
\begin{eqnarray}
\nonumber
 \|\overline{\widehat{\psi}}\widetilde{\psi}\|_{L^{1+}(I,H^{m-1,2})} & \le & c \|\widehat{\psi}\|_{L^{\infty}(I,H^{s,2})}\|\widetilde{\psi}\|_{L^{\infty}(I,L^2)}^{s-m+\frac{1}{2}} \|\widetilde{\psi}\|_{L^{\infty}(I,H^{1,2})}^{\frac{1}{2}+m-s} |I|^{1-} \\
\nonumber
& \le & c N^{(1-s\wedge m)(\frac{1}{2}+m-s)} N^{-\frac{4}{3}(1-s\wedge m)-\delta+} \\
\label{4.18}
& \le & c N^{(1-s\wedge m)(-\frac{5}{6}+m-s)+} \\
\nonumber
& \le & c N^{-\gamma (s,m)}
\end{eqnarray}
where $ \gamma (s,m) > 0 $.\\ 
In the same way we estimate $G_3$. Concerning $G_1$ we use (\ref{4.16}) with $ l=m-\frac{3}{2}$, $ \gamma = \frac{2}{\frac{3}{2}-m}$ , $ \rho = \frac{2}{m-\frac{1}{2}}$ and get by the embeddings $H^{s+,p} \subset B^{s,p} \subset H^{s-,p} $ (cf. \cite{T}, p. 180) and the definition of the Besov spaces $B^{s,p}$ (cf. \cite{BL}, Thm. 6.2.5) (alternatively one could also use the so-called fractional Leibniz rule):
\begin{eqnarray*}
 \lefteqn{ \| | \widehat{\psi}|^2\|_{H^{\frac{1}{2}+\epsilon,\frac{2}{\frac{5}{2}-m}+}} \le  c\| | \widehat{\psi}|^2\|_{B^{\frac{1}{2}+2\epsilon,\frac{2}{\frac{5}{2}-m}+}} } \\
& = : & c(\| |\widehat{\psi}|^2\|_{L^{\frac{2}{\frac{5}{2}-m}+}} + (\int_0^{\infty} (\tau^{-\frac{1}{2}-2\epsilon} \sup_{|h|\le \tau} \||\widehat{\psi}|^2(\cdot +h)-|\widehat{\psi}|^2(\cdot)\|_{L^{\frac{2}{\frac{5}{2}-m}+}})^2 \frac{d\tau}{\tau})^{\frac{1}{2}}) \\
& \le & c(\| \widehat{\psi}\|_{L^{\frac{2}{\frac{5}{2}-m}\widehat{p}+}} \| \widehat{\psi}\|_{L^{\frac{2}{\frac{5}{2}-m}\widehat{q}+}}  \\ 
& & + (\int_0^{\infty} (\tau^{-\frac{1}{2}-2\epsilon} \sup_{|h|\le \tau} \|\widehat{\psi}(\cdot +h)-\widehat{\psi}(\cdot)\|_{L^{\frac{2}{\frac{5}{2}-m}\widehat{q}+}})^2 \frac{d\tau}{\tau})^{\frac{1}{2}} \|\widehat{\psi}\|_{L^{\frac{2}{\frac{5}{2}-m}\widehat{p}+}}) \\
& \le & c\| \widehat{\psi}\|_{L^{\frac{2}{\frac{5}{2}-m}\widehat{p}+}}( \| \widehat{\psi}\|_{L^{\frac{2}{\frac{5}{2}-m}\widehat{q}+}} + \|\widehat{\psi}\|_{B^{\frac{1}{2}+2\epsilon,\frac{2}{\frac{5}{2}-m}\widehat{q}+}}) \\
& \le & c\| \widehat{\psi}\|_{L^{\frac{2}{\frac{5}{2}-m}\widehat{p}+}} \|\widehat{\psi}\|_{H^{\frac{1}{2}+3\epsilon,\frac{2}{\frac{5}{2}-m}\widehat{q}+}}
\end{eqnarray*}
This implies 
\begin{equation}
\label{N}
 \| | \widehat{\psi}|^2\|_{L^{\frac{2}{\frac{1}{2}+m}+}(I,H^{\frac{1}{2}+,\frac{2}{\frac{5}{2}-m}+})} \le c \| \widehat{\psi}\|_{L^{\infty}(I,L^{\frac{2}{\frac{5}{2}-m}\widehat{p}+})} \| \widehat{\psi}\|_{L^{\infty}(I,H^{\frac{1}{2}+,\frac{2}{\frac{5}{2}-m}\widehat{q}+})} |I|^{\frac{\frac{1}{2}+m}{2}-}
\end{equation}
The H\"older exponents $\widehat{p},\widehat{q}$ are chosen such that
 $$ H^{s,2} \subset L^{\frac{2}{\frac{5}{2}-m}\widehat{p}+} \cap H^{\frac{1}{2}+,\frac{2}{\frac{5}{2}-m}\widehat{q}+} $$
This requires 
$ \frac{1}{2} > \frac{1}{\widehat{p}} \cdot \frac{\frac{5}{2}-m}{2} > \frac{1}{2} - \frac{s}{3}$ and $ \frac{1}{2} > \frac{1}{\widehat{q}} \cdot \frac{\frac{5}{2}-m}{2} > \frac{1}{2} - \frac{s-\frac{1}{2}}{3} $
which can be fulfilled if
$ \frac{2}{\frac{5}{2}-m} > 1 = \frac{1}{\widehat{p}} + \frac{1}{\widehat{q}} > \frac{2}{\frac{5}{2}-m} \cdot (\frac{1}{2}-\frac{s}{3} + \frac{1}{2} - \frac{s-\frac{1}{2}}{3}) = \frac{2}{\frac{5}{2}-m} \cdot (\frac{7}{6} - \frac{2}{3}s) $.
The first inequality holds for $ m> \frac{1}{2} $ and the second one, because for $m\le 1$ , $ s > \frac{5}{8}$:
$ \frac{2}{\frac{5}{2}-m} \cdot (\frac{7}{6} - \frac{2}{3}s) < \frac{2}{\frac{3}{2}}(\frac{7}{6} - \frac{5}{12}) = \frac{4}{3} \cdot \frac{9}{12} = 1 $
Thus
\begin{eqnarray}
\nonumber
\| |\widehat{\psi}|^2\|_{L^{\frac{2}{\frac{1}{2}+m}+}(I,H^{\frac{1}{2}+,\frac{2}{\frac{5}{2}-m}+})} & \le & c \| \widehat{\psi} \|^2_{L^{\infty}(I,H^{s,2})} |I|^{\frac{\frac{1}{2}+m}{2}-} \\
\label{4.19}
& \le & c N^{-\frac{4}{3}(1-s\wedge m)\frac{\frac{1}{2}+m}{2}} \le c N^{-\gamma (s,m)}
\end{eqnarray}
From (\ref{4.2}),(\ref{4.17}),(\ref{4.18}),(\ref{4.19}) we conclude
\begin{eqnarray*}
\| S_{\pm}\widehat{\phi}_{\pm} \|_{X^{m,\frac{1}{2}+}(I)}  &
 \le & c \|\phi_{0\pm 2} \|_{H^{m,2}} + c \|G\|_{X^{m,-\frac{1}{2}+}_{\pm}(I)} \\
 & \le & c \|\phi_{0\pm 2} \|_{H^{m,2}} +c(\| \overline{\widehat{\psi}} \widetilde{\psi}\|_{L^{1+}(I,H^{m-1,2-})} + \| \widehat{\psi} \overline{\widetilde{\psi}}\|_{L^{1+}(I,H^{m-1,2-})} \\ & & \qquad \qquad \qquad \quad + \| |\widehat{\psi}|^2\|_{L^{\frac{2}{\frac{1}{2}+m}+}(I,H^{\frac{1}{2}+,\frac{2}{\frac{5}{2}-m}+})}) \\ 
& \le & c \|\phi_{0\pm 2} \|_{H^{m,2}} + c N^{-\gamma (s,m)} \\
& \le & c_0
\end{eqnarray*} 
provided $ c_0 \ge 2c \|\phi_{0\pm 2}\|_{H^{m,2}} $ and $N$ is sufficiently large.\\
Finally we treat $\| S_{\pm}\widehat{\phi}_{\pm} \|_{X_{\pm}^{0,\frac{1}{2}+}(I)}$. We use (\ref{3.5}) and estimate by interpolation
\begin{eqnarray}
\nonumber
\| |\widehat{\psi}|^2 \|_{L^{1+}(I,H^{-1,2})} & \le & c \| |\widehat{\psi}|^2 \|_{L^{1+}(I,L^{\frac{6}{5}})} = c \|\widehat{\psi} \|^2_{L^{2+}(I,L^{\frac{12}{5}})} \\
\nonumber
& \le & c \| \widehat{\psi} \|_{L^{\infty}(I,L^2)}^{2(1-\frac{1}{4s})} \| \widehat{\psi}\|_{L^{\infty}(I,H^{s,2})}^{\frac{1}{2s}} |I|^{1-} \\
\nonumber
&\le & c N^{-2(s\wedge m)(1-\frac{1}{4s})} N^{-\frac{4}{3}(1-s\wedge m)-\delta} \le c N^{-\frac{4}{3}-(s\wedge m)(\frac{2}{3}-\frac{1}{2s})-\delta} \\
\label{4.20}
& \le & c N^{-\frac{4}{3} - (\frac{2}{3} -1)-\delta} = c N^{-1-\delta}
\end{eqnarray}
because $s,m \le 1$ , $ s \ge \frac{1}{2}$.\\
Next we estimate by Lemma \ref{Lemma 3.1}:
\begin{eqnarray*}
\lefteqn{ \|\overline{\widehat{\psi}}\widetilde{\psi}\|_{L^{1+}(I,H^{-1,2})} } \\
& & \le c \|\widehat{\psi}\widetilde{\psi}\|_{L^{1+}(I,L^{\frac{6}{5}})} \le c \|\widehat{\psi}\|_{L^{\infty}(I,L^2)} \|\widetilde{\psi}\|_{L^{\infty}(I,L^3)}|I|^{1-} \\
& & \le c \|\widehat{\psi}\|_{L^{\infty}(I,L^2)} \|\widetilde{\psi}\|_{L^{\infty}(I,L^2)}^{\frac{1}{2}} \|\widetilde{\psi}\|_{L^{\infty}(I,H^{1,2})}^{\frac{1}{2}} |I|^{1-} \\
& & \le c \|\widehat{\psi}\|_{L^{\infty}(I,L^2)} N^{\frac{1}{2}(1-s\wedge m)} N^{-\frac{4}{3}(1-s\wedge m)-\delta} \\
& & \le c \|\widehat{\psi}\|_{L^{\infty}(I,L^2)} N^{-\frac{5}{6}(1-s\wedge m)-\delta}
\end{eqnarray*}
\begin{enumerate}
\item If $ s \le m $ we get
\begin{equation}
\label{4.21}
 \|\overline{\widehat{\psi}}\widetilde{\psi}\|_{L^{1+}(I,H^{-1,2})} \le c N^{-s-\frac{5}{6}(1-s)-\delta} \le c N^{-\frac{5}{6}-\frac{1}{6}s-\delta} 
\end{equation}
\item If $ m \le s $ we get
\begin{eqnarray}
\nonumber
\|\overline{\widehat{\psi}}\widetilde{\psi}\|_{L^{1+}(I,H^{-1,2})} & \le & c(N^{-\frac{1}{6}-\frac{5}{6}m} + N^{-s})N^{-\frac{5}{6}(1-m)-\delta} \\
\label{4.22}
& = & c(N^{-1-\delta} + N^{-s-\frac{5}{6}+\frac{5}{6}m-\delta}) \le c N^{-m-\delta}
\end{eqnarray}
because $-s-\frac{5}{6}+\frac{5}{6}m \le -m \Leftrightarrow \frac{11}{6}m \le s + \frac{5}{6} $. This is fulfilled, because $ \frac{11}{6}m \le \frac{11}{6}s = s + \frac{5}{6}s \le s + \frac{5}{6}$.
\end{enumerate}
The term $\|\widehat{\psi}\overline{\widetilde{\psi}}\|_{L^{1+}(I,H^{-1,2})}$ is estimated in the same way.\\
From (\ref{4.2}),(\ref{4.20}),(\ref{4.21}),(\ref{4.2}) we get 
\begin{enumerate}
\item in the case $s\le m$:
\begin{eqnarray*}
\|S_{\pm}\widehat{\phi}_{\pm}\|_{X^{0,\frac{1}{2}+}_{\pm}(I)} & \le & c \|\phi_{0\pm 2}\|_{L^2} + c N^{-\frac{5}{6}-\frac{1}{6}s-} \\
& \le & c N^{-m} + c N^{-\frac{5}{6}- \frac{1}{6}s-} \\
& \le & c_0 (N^{-m} + N^{-\frac{5}{6}-\frac{1}{6}s})
\end{eqnarray*}
if $ c_0 \ge c $ and $N$ sufficiently large.
\item Similarly, in the case $m\le s$:
$$ \|S_{\pm}\widehat{\phi}_{\pm}\|_{X^{0,\frac{1}{2}+}_{\pm}(I)} \le c_0 N^{-m} $$
\end{enumerate}
Summarizing, we have shown that $S$ maps $Z$ into itself. The contraction property uses exactly the same type of estimates and is therefore omitted.\\[2ex]
The next estimates show that the nonlinear part $w$ of (\ref{4.1}) behaves better than the linear part.
\begin{prop}
\label{Proposition 4.2}
Under the assumptions of Proposition \ref{Proposition 4.1} the following estimates hold:
\begin{eqnarray*}
\|w\|_{X^{1,\frac{1}{2}+}(I)} & \le & c N^{\frac{5}{6}-\frac{4}{3}(s\wedge m)} \\ 
\|w\|_{X^{0,\frac{1}{2}+}(I)} & \le & c N^{\frac{2}{3}-\frac{5}{3}(s\wedge m)}
\end{eqnarray*}
if in addition we assume $ s+m > \frac{3}{2} $.
\end{prop}
{\bf Proof:}
By Lemma \ref{Lemma 1.2} and Lemma \ref{Lemma 3.1} we have
\begin{eqnarray}
\label{4.22a}
\|w\|_{X^{0,\frac{1}{2}+}(I)} & \le & c \|F\|_{X^{0,-\frac{1}{4}-}(I)} |I|^{\frac{1}{4}-} \\
\nonumber
& \le & c(\|\widetilde{\phi}_{\pm} \widehat{\psi}\|_{X^{0,-\frac{1}{4}-}(I)} + \|\widehat{\phi}_{\pm} \widehat{\psi}\|_{X^{0,-\frac{1}{4}-}(I)} + \|\widetilde{\phi}_{\pm} \widetilde{\psi}\|_{X^{0,-\frac{1}{4}-}(I)}) |I|^{\frac{1}{4}-} \\
\nonumber
& \le & c(\|\widetilde{\phi}_{\pm}\|_{X^{0,\frac{1}{2}+}_{\pm}(I)} \| \widehat{\psi}\|_{X^{0,\frac{1}{2}+}(I)} + \|\widehat{\phi}_{\pm}\|_{X^{0,\frac{1}{2}+}_{\pm}(I)} \|\widehat{\psi}\|_{X^{0,\frac{1}{2}+}(I)} \\ \nonumber
& & \quad  + \|\widehat{\phi}_{\pm}\|_{X^{0,\frac{1}{2}+}_{\pm}(I)} \|\widetilde{\psi}\|_{X^{0,\frac{1}{2}+}(I)}) |I|^{\frac{1}{4}-} \\
\nonumber
& \le & c (N^{1-s\wedge m} N^{-s\wedge m} + N^{-s\wedge m} N^{-s\wedge m} + N^{-s\wedge m}) N^{-\frac{1}{3}(1-s\wedge m)-\frac{\delta}{4}-} \\
\nonumber
& \le & c N^{\frac{2}{3}-\frac{5}{3}(s\wedge m)}
\end{eqnarray}
In order to estimate $\|w\|_{X^{1,\frac{1}{2}+}(I)} $ we consider separately the cases $s\le m$ and $m\le s$.\\
Case $s\le m$:\\
By Lemma \ref{Lemma 0.5} with $s=1+$ and $\sigma=\frac{1}{2}-$ we get by interpolation
\begin{eqnarray}
\nonumber
\|\widetilde{\phi}_{\pm} \widehat{\psi}\|_{X^{1+,-\frac{1}{2}-}(I)} & \le & c( \|\widetilde{\phi}_{\pm}\|_{X_{\pm}^{\frac{1}{2}+,\frac{1}{2}+}(I)} \|\widehat{\psi}\|_{X^{0,\frac{1}{2}+}(I)} + \|\widetilde{\phi}_{\pm}\|_{X_{\pm}^{1,\frac{1}{2}+}(I)} \|\widehat{\psi}\|_{X^{\frac{1}{2}+,0}(I)}) \\
\nonumber
& \le & c( \|\widetilde{\phi}_{\pm}\|_{X_{\pm}^{\frac{1}{2}+,\frac{1}{2}+}(I)} \|\widehat{\psi}\|_{X^{0,\frac{1}{2}+}(I)} \\
\label{4.23}
& & \quad + \|\widetilde{\phi}_{\pm}\|_{X_{\pm}^{1,\frac{1}{2}+}(I)} \|\widehat{\psi}\|_{X^{0,\frac{1}{2}+}(I)}^{1-\frac{1}{2s}-}\|\widehat{\psi}\|_{X^{s,\frac{1}{2}+}(I)}^{\frac{1}{2s}+}|I|^{\frac{1}{2}})
\end{eqnarray}
because
\begin{eqnarray*}
\|\widehat{\psi}\|_{X^{\frac{1}{2}+,0}(I)} & = & \left(\int_I \|\widehat{\psi}(t)\|_{H^{\frac{1}{2}+,2}}^2 dt\right)^{\frac{1}{2}} \le \|\widehat{\psi}\|_{L^{\infty}(I,H^{\frac{1}{2}+,2})} |I|^{\frac{1}{2}} \\
& \le & c \|\widehat{\psi}\|_{X^{\frac{1}{2}+,\frac{1}{2}+}(I)} |I|^{\frac{1}{2}} \le c \|\widehat{\psi}\|_{X^{0,\frac{1}{2}+}(I)}^{1-\frac{1}{2s}-} \|\widehat{\psi}\|_{X^{s,\frac{1}{2}+}(I)}^{\frac{1}{2s}+} |I|^{\frac{1}{2}}
\end{eqnarray*}
Thus
\begin{eqnarray*}
\|\widetilde{\phi}_{\pm}\widehat{\psi}\|_{X^{1+,-\frac{1}{2}-}(I)} & \le & c(N^{1-s\wedge m} N^{-s} + N^{1-s\wedge m} N^{-s(1-\frac{1}{2s})+} N^{-\frac{2}{3}(1-s\wedge m)-\frac{\delta}{2}+}) \\
& \le & c(N^{1-2s} + N^{1-s-s+\frac{1}{2}-\frac{2}{3}+\frac{2}{3}s-\frac{\delta}{2}+}) \\
& \le & c N^{\frac{5}{6}-\frac{4}{3}s-\frac{\delta}{2}+}
\end{eqnarray*}
Next from Lemma \ref{Lemma 0.5} with $s=1+$ and $\sigma=0+$, Lemma \ref{Lemma 3.1} and interpolation
\begin{eqnarray}
\nonumber
\lefteqn{\|\widehat{\phi}_{\pm}\widetilde{\psi}\|_{X^{1+,-\frac{1}{2}-}(I)}} \\ \nonumber
& & \le  c(\|\widetilde{\psi}\|_{X^{0,\frac{1}{2}+}(I)} \|\widehat{\phi}_{\pm}\|_{X_{\pm}^{\frac{1}{2}+,\frac{1}{2}+}(I)} + \|\widetilde{\psi}\|_{X^{1,0}(I)} \|\widehat{\phi}_{\pm}\|_{X_{\pm}^{\frac{1}{2}+,\frac{1}{2}+}(I)}) \\
\nonumber
& & \le  c(\|\widetilde{\psi}\|_{X^{0,\frac{1}{2}+}(I)} \|\widehat{\phi}_{\pm}\|_{X_{\pm}^{0,\frac{1}{2}+}(I)}^{1-\frac{1}{2m}-} \|\widehat{\phi}_{\pm}\|_{X_{\pm}^{m,\frac{1}{2}+}(I)}^{\frac{1}{2m}+} \\ \nonumber
& & \quad \quad + \|\widetilde{\psi}\|_{X^{1,\frac{1}{2}+}(I)} |I|^{\frac{1}{2}} \|\widehat{\phi}_{\pm}\|_{X_{\pm}^{0,\frac{1}{2}+}(I)}^{1-\frac{1}{2m}-}\|\widehat{\phi}_{\pm}\|_{X_{\pm}^{m,\frac{1}{2}+}(I)}^{\frac{1}{2m}+}) \\
\label{4.24}
& & \le  c(N^{-(s\wedge m)(1-\frac{1}{2m})+} + N^{1-s\wedge m} N^{-\frac{2}{3}(1-s\wedge m)-\frac{\delta}{2}+} N^{-(s\wedge m)(1-\frac{1}{2m})+}) \\
\nonumber
& & \le  c(N^{-(s\wedge m)(1-\frac{1}{2m})+} + N^{\frac{1}{3}-(s\wedge m)(\frac{4}{3}-\frac{1}{2m})-\frac{\delta}{2}+}) \\
\nonumber
& & \le  c N^{\frac{1}{3}-(s\wedge m)(\frac{4}{3}-\frac{1}{2m})-\frac{\delta}{2}+} \\
\nonumber
& & \le c N^{\frac{1}{3}-(s\wedge m)(\frac{4}{3}-\frac{1}{2(s\wedge m)})-\frac{\delta}{2}+} \\
\nonumber
& & = cN^{\frac{5}{6}-\frac{4}{3}(s\wedge m)-\frac{\delta}{2}+} \\
\nonumber
& & = cN^{\frac{5}{6}-\frac{4}{3}s-\frac{\delta}{2}+}
\end{eqnarray}
because $-(s\wedge m)(1-\frac{1}{2m})<\frac{1}{3}-(s\wedge m)(\frac{4}{3}-\frac{1}{2m}) \Leftrightarrow (s\wedge m)\frac{1}{3} < \frac{1}{3}$ which is fulfilled.\\
Moreover by Lemma \ref{Lemma 0.5} with $s=1+$ and $\sigma=m-\frac{1}{2}-$ and interpolation
\begin{eqnarray*}
\lefteqn{ \|\widehat{\phi}_{\pm}\widehat{\psi}\|_{X^{1+,-\frac{1}{2}-}(I)} } \\
& & \le c(\|\widehat{\psi}\|_{X^{0,\frac{1}{2}+}(I)} \|\widehat{\phi}_{\pm}\|_{X_{\pm}^{0,\frac{1}{2}+}(I)}^{1-\frac{1}{2m}-} \|\widehat{\phi}_{\pm}\|_{X_{\pm}^{m,\frac{1}{2}+}(I)}^{\frac{1}{2m}+} + \|\widehat{\psi}\|_{X^{\frac{3}{2}-m+,0}(I)} \|\widehat{\phi}_{\pm}\|_{X_{\pm}^{m,\frac{1}{2}+}(I)})
\end{eqnarray*}
using
\begin{eqnarray*}
\|\widehat{\psi}\|_{X^{\frac{3}{2}-m+,0}(I)} & = & \|\widehat{\psi}\|_{L^2_t(I,H^{\frac{3}{2}-m+,2}_x)} \le \|\widehat{\psi}\|_{L^{\infty}_t(I,H^{\frac{3}{2}-m+,2}_x)}|I|^{\frac{1}{2}} \\
& \le & c \|\widehat{\psi}\|_{X^{\frac{3}{2}-m+,\frac{1}{2}+}(I)} |I|^{\frac{1}{2}} \le c \|\widehat{\psi}\|_{X^{0,\frac{1}{2}+}(I)}^{1-\frac{\frac{3}{2}-m}{s}-} \|\widehat{\psi}\|_{X^{s,\frac{1}{2}+}(I)}^{\frac{\frac{3}{2}-m}{s}+} |I|^{\frac{1}{2}}
\end{eqnarray*}
under our assumption $ s+m > \frac{3}{2} $.\\
Thus we get
\begin{eqnarray*}
\lefteqn{ \|\widehat{\phi}_{\pm}\widehat{\psi}\|_{X^{1+,-\frac{1}{2}-}(I)} } \\
& & \le c(N^{-s} N^{-(s\wedge m)(1-\frac{1}{2m})+} + N^{-s(1-\frac{\frac{3}{2}-m}{s})+} N^{-\frac{2}{3}(1-s\wedge m)-\frac{\delta}{2}}) \\ 
& & \le c(N^{-s}N^{-(s\wedge m)(1-\frac{1}{2(s\wedge m)})+} + N^{-s+\frac{3}{2}-m-\frac{2}{3}+\frac{2}{3}s-\frac{\delta}{2}+})  \\
& & \le c(N^{\frac{1}{2}-2s+} + N^{\frac{5}{6}-\frac{4}{3}s-\frac{\delta}{2}+}) \\
& & \le c N^{\frac{5}{6}-\frac{4}{3}s-\frac{\delta}{2}+}
\end{eqnarray*}
Case $m\le s$:\\
From (\ref{4.23}) we conclude
\begin{eqnarray*}
\lefteqn{ \|\widetilde{\phi}_{\pm} \widehat{\psi}\|_{X^{1+,-\frac{1}{2}-}(I)} } \\
& & \le c (N^{1-s\wedge m} N^{-s\wedge m} + N^{1-s\wedge m} N^{-s\wedge m(1-\frac{1}{2s})+} N^{-\frac{2}{3}(1-s\wedge m)-\frac{\delta}{2}+}) \\
& & \le c (N^{1-2(s\wedge m)} + N^{\frac{1}{3}-(s\wedge m)(1+1-\frac{1}{2s}-\frac{2}{3})-\frac{\delta}{2}+}) \\
& & = c (N^{1-2(s\wedge m)} + N^{\frac{1}{3}-(s\wedge m)(\frac{4}{3}-\frac{1}{2s})-\frac{\delta}{2}+}) \\ 
& & \le c (N^{1-2(s\wedge m)} + N^{\frac{1}{3}-(s\wedge m)(\frac{4}{3}-\frac{1}{2(s\wedge m)})-\frac{\delta}{2}+}) \\
& & = c (N^{1-2(s\wedge m)} + N^{\frac{5}{6}-\frac{4}{3}(s\wedge m)-\frac{\delta}{2}+}) \\
& & \le c N^{\frac{5}{6}-\frac{4}{3}m-\frac{\delta}{2}+}
\end{eqnarray*}
because $1-2(s\wedge m) < \frac{5}{6}-\frac{4}{3}(s\wedge m) \Leftrightarrow 1<4(s\wedge m)$ which is fulfilled.\\
Next from (\ref{4.24}):
\begin{eqnarray*}
\lefteqn{\|\widehat{\phi}_{\pm}\widetilde{\psi}\|_{X^{1+,-\frac{1}{2}-}(I)}} \\
& & \le c (N^{-m(1-\frac{1}{2m})+} + N^{1-s\wedge m} N^{-\frac{2}{3}(1-s\wedge m)-\frac{\delta}{2}+} N^{-m(1-\frac{1}{2m}+}) \\
& & \le c(N^{-m+\frac{1}{2}+} + N^{1-m-\frac{2}{3}+\frac{2}{3}m-m+\frac{1}{2}-\frac{\delta}{2}+}) \\
& & \le c N^{\frac{5}{6}-\frac{4}{3}m-\frac{\delta}{2}+}
\end{eqnarray*}
because $-m+\frac{1}{2} < \frac{5}{6} - \frac{4}{3}m \Leftrightarrow m<1$.\\
Next by Lemma \ref{Lemma 0.5} with $s=1+$ , $ \sigma = 1-s+$ and interpolation, again using $s+m > \frac{3}{2}$:
\begin{eqnarray*}
\lefteqn{ \|\widehat{\phi}_{\pm} \widehat{\psi}\|_{X^{1+,-\frac{1}{2}-}(I)} } \\
& & \le c(\|\widehat{\psi}\|_{X^{0,\frac{1}{2}+}(I)} \|\widehat{\phi}_{\pm}\|_{X_{\pm}^{0,\frac{1}{2}+}(I)}^{1-\frac{1}{2m}-}
\|\widehat{\phi}_{\pm}\|_{X_{\pm}^{m,\frac{1}{2}+}(I)}^{\frac{1}{2m}+} + \|\widehat{\psi}\|_{X^{s,0}(I)} \|\widehat{\phi}_{\pm}\|_{X_{\pm}^{\frac{3}{2}-s+,\frac{1}{2}+}(I)}) \\
& & \le c(\|\widehat{\psi}\|_{X^{0,\frac{1}{2}+}(I)} \|\widehat{\phi}_{\pm}\|_{X_{\pm}^{0,\frac{1}{2}+}(I)}^{1-\frac{1}{2m}-}
\|\widehat{\phi}_{\pm}\|_{X_{\pm}^{m,\frac{1}{2}+}(I)}^{\frac{1}{2m}+} + \|\widehat{\psi}\|_{X^{s,\frac{1}{2}+}(I)}|I|^{\frac{1}{2}} \|\widehat{\phi}_{\pm}\|_{X_{\pm}^{\frac{3}{2}-s+,\frac{1}{2}+}(I)}) \\
& & \le c(\|\widehat{\psi}\|_{X^{0,\frac{1}{2}+}(I)} \|\widehat{\phi}_{\pm}\|_{X_{\pm}^{0,\frac{1}{2}+}(I)}^{1-\frac{1}{2m}-}
\|\widehat{\phi}_{\pm}\|_{X_{\pm}^{m,\frac{1}{2}+}(I)}^{\frac{1}{2m}+} \\
& & \qquad + \|\widehat{\psi}\|_{X^{s,\frac{1}{2}+}(I)} |I|^{\frac{1}{2}} \|\widehat{\phi}_{\pm}\|_{X_{\pm}^{0,\frac{1}{2}+}(I)}^{1-\frac{\frac{3}{2}-s}{m}-}\|\widehat{\phi}_{\pm}\|_{X_{\pm}^{m,\frac{1}{2}+}(I)}^{\frac{\frac{3}{2}-s}{m}+}) \\
& & \le c(N^{-s\wedge m} N^{-m(1-\frac{1}{2m})+} + N^{-\frac{2}{3}(1-s\wedge m)-\frac{\delta}{2}+} N^{-m(1-\frac{\frac{3}{2}-s}{m})+}) \\
& & \le c(N^{\frac{1}{2}-2m+} + N^{-m+\frac{3}{2}-s-\frac{2}{3}+\frac{2}{3}m-\frac{\delta}{2}+}) \\
& & \le c(N^{\frac{1}{2}-2m+} + N^{\frac{5}{6}-\frac{4}{3}m-\frac{\delta}{2}+}) \\
& & \le c N^{\frac{5}{6}-\frac{4}{3}m-\frac{\delta}{2}+}
\end{eqnarray*}
Summarizing we arrive at
\begin{equation}
\label{4.25}
\|F\|_{X^{1+,-\frac{1}{2}-}(I)} \le c N^{\frac{5}{6}-\frac{4}{3}(s\wedge m)-\frac{\delta}{2}+}
\end{equation}
Now by (\ref{4.22a}):
\begin{equation}
\label{4.26}
\|F\|_{X^{0,-\frac{1}{4}-}(I)} \le c N^{1-2(s\wedge m)}
\end{equation}
Interpolating between (\ref{4.25}) and (\ref{4.26}) we arrive at
$$ \|w\|_{X^{1,\frac{1}{2}+}(I)} \le c \|F\|_{X^{1,-\frac{1}{2}+}(I)} \le c N^{\frac{5}{6}-\frac{4}{3}(s\wedge m)} $$
The proof is complete.

The next result shows that the nonlinear part $z_{\pm}$ of (\ref{4.2}) also behaves better than the linear part.
\begin{prop}
\label{Proposition 4.3}
Under the assumptions of Proposition \ref{Proposition 4.1} the following estimates hold
\begin{eqnarray}
\label{4.27}
\|z_{\pm}\|_{L^{\infty}(I,H^{1,2})} & \le & c N^{-\frac{1}{3}-(s\wedge m)(\frac{2}{3}-\frac{1}{2s})}\\
\label{4.28}
\|z_{\pm}\|_{X_{\pm}^{0,\frac{1}{2}+}(I)} & \le & c N^{-(s\wedge m)}
\end{eqnarray}
\end{prop}
{\bf Proof:}
The estimate (\ref{4.28}) is already proven in Proposition \ref{Proposition 4.1}.\\
In order to prove (\ref{4.27}) we estimate by interpolation and Lemma \ref{Lemma 3.1}:
\begin{eqnarray}
\nonumber
\|\overline{\widehat{\psi}}\widetilde{\psi}\|_{L^1(I,L^2)} & \le & \|\widehat{\psi}\|_{L^{\infty}(I,L^3)} \|\widetilde{\psi}\|_{L^{\infty}(I,L^6)} |I| \\
\nonumber
& \le & c \|\widehat{\psi}\|_{L^{\infty}(I,L^2)}^{1-\frac{1}{2s}} \|\widehat{\psi}\|_{L^{\infty}(I,H^{s,2})}^{\frac{1}{2s}} \|\widetilde{\psi}\|_{L^{\infty}(I,H^{1,2})} |I| \\
\label{4.29}
& \le & c N^{-(s\wedge m)(1-\frac{1}{2s})} N^{1-s\wedge m} N^{-\frac{4}{3}(1-s\wedge m)-\delta} \\
\nonumber
& \le & c N^{-\frac{1}{3}-(s\wedge m)(\frac{2}{3}-\frac{1}{2s})-\delta}
\end{eqnarray}
The term $\|\widehat{\psi}\overline{\widetilde{\psi}}\|_{L^1(I,L^2)}$ is estimated in exactly the same way.\\
Finally we use Strichartz' inequality for the inhomogeneous Klein-Gordon equation $u_{tt} + Au = f$ , $ u(0)=u_t(0)=0 $ (cf. \cite{GV}):
$$ \|u\|_{L^{\infty}(I,H^{1,2}({\bf R}^3))} \le c \|f\|_{L^{\gamma '}(I,H^{-l,\rho '}({\bf R}^3))} $$
which holds if $ l+2(\frac{1}{2}-\frac{1}{\rho})=0$ and $\frac{1}{\gamma}=\frac{1}{2}-\frac{1}{\rho}$. Choose $l=-s+$ , $\rho=\frac{2}{1-s}-$ , $\gamma = \frac{2}{s}-$ , $\rho ' = \frac{2}{1+s}+$ , $\gamma ' = \frac{2}{2-s}+$ and estimate using Besov norms as in (\ref{N}) or alternatively by the fractional Leibniz rule:
$$ \| |\widehat{\psi}|^2 \|_{L^{\frac{2}{2-s}+}(I,H^{s-,\frac{2}{1+s}+})} \le c 
\| \widehat{\psi} \|_{L^{\infty}(I,L^{\frac{2}{1+s}\widehat{p}+})}           \|\widehat{\psi} \|_{L^{\infty}(I,H^{s,\frac{2}{1+s}\widehat{q}+})} |I|^{\frac{2-s}{2}-} $$
Choose $\widehat{p} = \frac{1+s}{s}+$ , $\widehat{q} = 1+s-$ and interpolate (for $ s> \frac{3}{5}$):
\begin{eqnarray*}
\| |\widehat{\psi}|^2 \|_{L^{\frac{2}{2-s}+}(I,H^{s-,\frac{2}{1+s}+})} & \le & c 
\| \widehat{\psi} \|_{L^{\infty}(I,L^{\frac{2}{s}+})}           \|\widehat{\psi} \|_{L^{\infty}(I,H^{s,2})} |I|^{\frac{2-s}{2}-} \\
& \le & c 
\| \widehat{\psi} \|_{L^{\infty}(I,L^2)}^{\frac{5s-3}{2s}-}           \|\widehat{\psi} \|_{L^{\infty}(I,H^{s,2})}^{2-\frac{5s-3}{2s}+} |I|^{\frac{2-s}{2}-} \\ 
& \le & c N^{-(s\wedge m)\frac{5s-3}{2s}+} N^{-\frac{2-s}{2}\frac{4}{3}(1-s\wedge m)-}
\end{eqnarray*}
Now we have
$-(s\wedge m)\frac{5s-3}{2s} - \frac{2-s}{2} \cdot \frac{4}{3} (1-s\wedge m) < - \frac{1}{3} - (s\wedge m)(\frac{2}{3} - \frac{1}{2s}) $
if 
$ -(s\wedge m)\frac{5s-3}{2s} - \frac{2}{3} (1-s\wedge m)  <  - \frac{1}{3} - (s\wedge m)(\frac{2}{3} - \frac{1}{2s}) \Leftrightarrow (s\wedge m)(\frac{1}{s}-\frac{7}{6}) < \frac{1}{3} $.
This holds if
$ 1-\frac{7}{6}(s\wedge m) = (s\wedge m)(\frac{1}{s\wedge m} - \frac{7}{6}) < \frac{1}{3} \Leftrightarrow \frac{4}{7} < s\wedge m$.
This is fulfilled, because $\frac{7}{10} > \frac{4}{7} $. Thus the decisive estimate is (\ref{4.29}), which directly leads to our claim (\ref{4.27}).
\section{The iteration process}
In the preceding sections we constructed a solution of the problem (\ref{-1.1}),(\ref{-1.2}) with data (\ref{-1.3}) $(\psi_0,\phi_0,\phi_1)$ in the time interval $I=[0,|I|]$ with $|I|=N^{-\frac{4}{3}(1-s\wedge m)-}$. Namely, if we define $\psi:=\widehat{\psi}+\widetilde{\psi}$ , $\phi_{\pm}:=\widehat{\phi}_{\pm}+\widetilde{\phi}_{\pm}$ we see that $(\psi,\phi_+,\phi_-)$ solves the system (\ref{1.4}),(\ref{1.5}) with initial conditions $\psi(0)=\psi_0$ , $ \psi_{\pm}(0)=\phi_{0\pm}$. This problem is equivalent to the original system (\ref{-1.1}),(\ref{-1.2}),(\ref{-1.3}).The initial data are transformed by $\phi_{0\pm}=\phi_0\pm iA^{-1/2}\phi_1$ or, conversely, by $\phi_0=\frac{1}{2}(\phi_{0+}+\phi_{0-})$ , $\phi_1=-\frac{i}{2}A^{1/2}(\phi_{0+}-\phi_{0-})$. In order to continue the solution of (\ref{1.4}),(\ref{1.5}) we take as new initial data $(\widetilde{\psi}(|I|)+w(|I|),\widetilde{\phi}_+(|I|)+z_+(|I|),\widetilde{\phi}_-(|I|) + z_-(|I|))$ instead of $(\psi_{01},\phi_{0+1},\phi_{0-1})$. When we have shown that this problem has a solution $(\widetilde{\widetilde{\psi}},\widetilde{\widetilde{\phi}}_+,\widetilde{\widetilde{\phi}}_-)$ in the interval $[|I|,2|I|]$ with equal length $|I|$ we insert this solution into the system (\ref{4.0a}),(\ref{4.0b}) in place of $(\widetilde{\psi},\widetilde{\phi}_+,\widetilde{\phi}_-)$ and solve this problem with data $(e^{i|I|\Delta}\psi_{02},e^{-i|I|A^{1/2}}\phi_{0+2},e^{i|I|A^{1/2}}\phi_{0-2})$ in $[|I|,2|I|]$.The solution of the original system (\ref{-1.1}),(\ref{-1.2}) corresponding to $(\widetilde{\widetilde{\psi}},\widetilde{\widetilde{\phi}}_+,\widetilde{\widetilde{\phi}}_-)$, denoted by $(\widetilde{\widetilde{\psi}},\widetilde{\widetilde{\phi}})$, then obviously has the following initial data (using (\ref{1.3a}),(\ref{1.3b})):
\begin{eqnarray*}
\widetilde{\widetilde{\psi}}(|I|) & = & \widetilde{\psi}(|I|) + w(|I|) \\
\widetilde{\widetilde{\phi}}(|I|) & = & \frac{1}{2}(\widetilde{\widetilde{\phi}}_+(|I|) + \widetilde{\widetilde{\phi}}_-(|I|)) = \frac{1}{2}(\widetilde{\phi}_+(|I|)+z_+(|I|)+\widetilde{\phi}_-(|I|)+z_-(|I|)) \\
& = & \widetilde{\phi}(|I|) + \frac{1}{2}(z_+(|I|)+z_-(|I|)) =: \widetilde{\phi}(|I|) + z(|I|) \\
\widetilde{\widetilde{\phi_t}}(|I|) & = & -\frac{i}{2}A^{\frac{1}{2}}(\widetilde{\widetilde{\phi}}_+(|I|)-\widetilde{\widetilde{\phi}}_-(|I|)) \\
& = & -\frac{i}{2}A^{\frac{1}{2}}(\widetilde{\phi}_+(|I|)+z_+(|I|)-\widetilde{\phi}_-(|I|)-z_-(|I|)) \\
& = & \widetilde{\phi}_t(|I|) - \frac{i}{2}A^{\frac{1}{2}}(z_+(|I|)-z_-(|I|)) =: \widetilde{\phi}_t(|I|) + z'(|I|)
\end{eqnarray*}
Adding up the solutions, we get a solution of the original problem in $[|I|,2|I|]$ as before. This defines an iteration process. At each step we have to ensure the same bounds on the initial data which were used in the first step. The replacement of $(\psi_{02},\phi_{0+2},\phi_{0-2})$ by $(e^{i|I|\Delta}\psi_{02},e^{-i|I|A^{1/2}}\phi_{0+2},e^{i|I|A^{1/2}}\phi_{0-2})$ is harmless, because the $H^s$-norms remain unchanged. The bounds on the data are controlled by the energy and the $L^2$-conservation law (see (\ref{2.1}),(\ref{2.5}),(\ref{2.6})). Thus we have to estimate these quantities independently of the iteration step. This is easy for the $L^2$-conserved quantity, the increment when replacing $\psi_{01}$ by $\widetilde{\psi}(|I|)+w(|I|)$ is given by
\begin{eqnarray*}
\left| \|\widetilde{\psi}(|I|)+w(|I|)\|_{L^2_x} - \|\psi_{01}\|_{L^2_x} \right| & =	& \left| \|\widetilde{\psi}(|I|)+w(|I|)\|_{L^2_x} - \|\widetilde{\psi}(|I|)\|_{L^2_x} \right| \\
& \le & \|w(|I|)\|_{L^2_x} \le c_2 N^{\frac{2}{3}-\frac{5}{3}(s\wedge m)} 
\end{eqnarray*}
by Proposition \ref{Proposition 4.2}, where $c_2=c_2(\overline{c},M)$.\\
The number of iteration steps in order to reach the given time $T$ is $\frac{T}{|I|}=TN^{\frac{4}{3}(1-s\wedge m)+}$. This means that in order to get uniform control over the $L^2$-norm of $\widetilde{\psi},\widetilde{\widetilde{\psi}},...$ we have to ensure that $c_2TN^{\frac{4}{3}(1-s\wedge m)+}N^{\frac{2}{3}-\frac{5}{3}(s\wedge m)}<M$, where $c_2=c_2(2\overline{c},2M)$ (remark that initially the $L^2$-norm of $\widetilde{\psi}$ was bounded by $M$). This is fulfilled for $N$ sufficiently large if $\frac{4}{3}(1-s\wedge m)+\frac{2}{3}-\frac{5}{3}(s\wedge m)<0 \Leftrightarrow s\wedge m > \frac{2}{3}$, which is fulfilled.\\
Concerning the increment of the energy, we define $z=\frac{1}{2}(z_++z_-)$ and estimate:
\begin{eqnarray*}
& & \hspace{-0.5cm} |E(\widetilde{\psi}(|I|)+w(|I|),\widetilde{\phi}(|I|)+z(|I|),\widetilde{\phi}_t(|I|)+z'(|I|)) - E(\psi_{01},\phi_{01},\phi_{11})| \\
& \hspace{-0.5cm} = & \hspace{-0.2cm} |E(\widetilde{\psi}(|I|)+w(|I|),\widetilde{\phi}(|I|)+z(|I|),\widetilde{\phi}_t(|I|)+z'(|I|)) - E(\widetilde{\psi}(|I|),\widetilde{\phi}(|I|),\widetilde{\phi}_t(|I|))| \\
& \hspace{-0.5cm} \le & \hspace{-0.2cm} 2(\|\nabla \widetilde{\psi}(|I|)\|+\|\nabla w(|I|)\|)\|\nabla w(|I|)\| + (\|A^{\frac{1}{2}}\widetilde{\phi}(|I|)\|+\|A^{\frac{1}{2}}z(|I|)\|)\|A^{\frac{1}{2}}z(|I|)\| \\
& & \quad + (\|\widetilde{\phi}_t(|I|)\|+\|z'(|I|)\|)\|z'(|I|)\| + \int_{{\bf R}^3} |z(|I|)| \, |\widetilde{\psi}(|I|)+w(|I|)|^2 dx \\
& & \quad + \int_{{\bf R}^3} |\widetilde{\phi}(|I|)| \, \left| |\widetilde{\psi}(|I|)+w(|I|)|^2 - |\widetilde{\psi}(|I|)|^2\right| dx
\end{eqnarray*}
The first term is bounded by Proposition \ref{Proposition 4.2} and (\ref{2.10}) by
$$c(N^{1-s\wedge m} + N^{\frac{5}{6}-\frac{4}{3}(s\wedge m)})N^{\frac{5}{6}-\frac{4}{3}(s\wedge m)} \le c N^{1-s\wedge m} N^{\frac{5}{6}-\frac{4}{3}(s\wedge m)} = cN^{\frac{11}{6}-\frac{7}{3}(s\wedge m)}$$
The second and third term are bounded by
$$c(N^{1-s\wedge m} + N^{-\frac{1}{3}-(s\wedge m)(\frac{2}{3}-\frac{1}{2s})}) N^{-\frac{1}{3}-(s\wedge m)(\frac{2}{3}-\frac{1}{2s})}$$
using (\ref{2.9}) and (\ref{4.27}). Now $1-s\wedge m \ge -\frac{1}{3}-(s\wedge m)(\frac{2}{3}-\frac{1}{2s}) \Leftrightarrow \frac{4}{3} \ge (s\wedge m)(\frac{1}{3}+\frac{1}{2s})$ which is fulfilled for $1\ge s,m \ge \frac{1}{2}$.
Consequently the following bound holds for the second and third term: $cN^{\frac{2}{3}-(s\wedge m)(\frac{5}{3}-\frac{1}{2s})}$.\\
Using (\ref{4.28}) and Proposition \ref{Proposition 4.2} the forth term is estimated as follows:
\begin{eqnarray*}
\int_{{\bf R}^3} |z(|I|)| | \widetilde{\psi}(|I|)+w(|I|)|^2 dx & \le & c \|z(|I|)\|_{L^2}(\|\widetilde{\psi}(|I|)\|_{L^4}^2 + \|w(|I|)\|_{L^4}^2) \\
& \le & c \|z(|I|)\|_{L^2}(\|\widetilde{\psi}(|I|)\|_{H^{1,2}}^2 + \|w(|I|)\|_{H^{1,2}}^2) \\
& \le & c N^{-(s\wedge m)}(N^{2(1-s\wedge m)} + N^{2(\frac{5}{6}-\frac{4}{3}(s\wedge m))}) \\
& \le & c N^{-(s\wedge m)}N^{2(1-s\wedge m)} \\
& \le & c N^{2-3(s\wedge m)}
\end{eqnarray*}
Finally the fifth term is estimated by (\ref{2.9}), Prop. \ref{Proposition 4.2} and (\ref{2.8}) as follows:
\begin{eqnarray*}
& & \int_{{\bf R}^3} |\widetilde{\phi}(|I|)| \, \left| |\widetilde{\psi}(|I|)+w(|I|)|^2 - |\widetilde{\psi}(|I|)|^2\right| dx \\
& \le & \|\widetilde{\phi}(|I|)\|_{L^6}\|w(|I|)\|_{L^2}(\|\widetilde{\psi}(|I|)\|_{L^3} + \|w(|I|)\|_{L^3}) \\
& \le & \|\nabla \widetilde{\psi}(|I|)\|_{L^2} \|w(|I|)\|_{L^2}(\|\widetilde{\psi}(|I|)\|_{L^2}^{\frac{1}{2}}\|\nabla \widetilde{\psi}(|I|)\|_{L^2}^{\frac{1}{2}} + \|w(|I|)\|_{L^2}^{\frac{1}{2}}\|\nabla w(|I|)\|_{L^2}^{\frac{1}{2}}) \\
& \le & c N^{1-s\wedge m} N^{\frac{2}{3}-\frac{5}{3}(s\wedge m)}(N^{\frac{1}{2}(1-s\wedge m)} + N^{\frac{1}{2}(\frac{2}{3}-\frac{5}{3}(s\wedge m))} N^{\frac{1}{2}(\frac{5}{6}-\frac{4}{3}(s\wedge m))}) \\ 
& = & c N^{1-s\wedge m} N^{\frac{2}{3}-\frac{5}{3}(s\wedge m)}(N^{\frac{1}{2}-\frac{1}{2}(s\wedge m)} + N^{\frac{3}{4}-\frac{3}{2}(s\wedge m)}) \\
& \le & c N^{1-s\wedge m} N^{\frac{2}{3}-\frac{5}{3}(s\wedge m)} N^{\frac{1}{2}-\frac{1}{2}(s\wedge m)} \\
& = & c N^{\frac{13}{6}-\frac{19}{6}(s\wedge m)}
\end{eqnarray*}
because $ \frac{1}{2}-\frac{1}{2}(s\wedge m) > \frac{3}{4}-\frac{3}{4}(s\wedge m)$ if $s\wedge m > \frac{1}{4}$.\\
Now, the forth term behaves better than the first one, because \\ 
$2-3(s\wedge m) < \frac{11}{6} - \frac{7}{3}(s\wedge m) \Leftrightarrow \frac{1}{6} < \frac{2}{3}(s\wedge m)$ which holds for $s\wedge m > \frac{1}{4}$.\\
Similarly, the fifth term is harmless compared to the first one, because \\
 $ \frac{13}{6}-\frac{19}{6}(s\wedge m) < \frac{11}{6} - \frac{7}{3}(s\wedge m) \Leftrightarrow \frac{2}{6} < \frac{5}{6}(s\wedge m) \Leftrightarrow s\wedge m > \frac{2}{5} $.\\
Thus the decisive terms are the first, second and third one.\\
Concerning the first term the condition that ensures uniform control of the energy of $(\widetilde{\psi},\widetilde{\phi}),(\widetilde{\widetilde{\psi}},\widetilde{\widetilde{\phi}}),...$ is the following:
$$ c_3TN^{\frac{4}{3}(1-s\wedge m)+} N^{\frac{11}{6}-\frac{7}{3}(s\wedge m)} < \overline{c} N^{2(1-s\wedge m)} $$
where $ c_3=c_3(2\overline{c},2M) $ (recall that the energy initially is bounded by $\overline{c}N^{2(1-s\wedge m)}$).\\
This is satisfied for $N$ sufficiently large provided \\ 
$\frac{4}{3}(1-s\wedge m)+\frac{11}{6}-\frac{7}{3}(s\wedge m) < 2(1-s\wedge m) \Leftrightarrow \frac{7}{6}<\frac{5}{3}(s\wedge m) \Leftrightarrow s\wedge m > \frac{7}{10} $.
Here is the point where the decisive bound on $s\wedge m$ appears.\\
Concerning the second and third term the following condition has to be satisfied:
$$ c_3TN^{\frac{4}{3}(1-s\wedge m)+} N^{\frac{2}{3}-(s\wedge m)(\frac{5}{3}-\frac{1}{2s})} < \overline{c} N^{2(1-s\wedge m)} $$
This requires\\
$ \frac{4}{3}(1-s\wedge m)+\frac{2}{3}-(s\wedge m)(\frac{5}{3}-\frac{1}{2s}) < 2(1-s\wedge m) \Leftrightarrow 0 < (s\wedge m)(1-\frac{1}{2s}) \Leftrightarrow s > \frac{1}{2} $.\\
The uniform control of the energy implies by (\ref{2.5}),(\ref{2.6}) uniform control of the $L^2$-norm of $(\nabla \widetilde{\psi},A^{\frac{1}{2}}\widetilde{\phi},\widetilde{\phi}_t)$,$(\nabla \widetilde{\widetilde{\psi}},A^{\frac{1}{2}}\widetilde{\widetilde{\phi}},\widetilde{\widetilde{\phi}}_t),...$ .\\

We have proven
\begin{theorem}
\label{Theorem 5.1}
Let $ 1 \ge s,m > \frac{7}{10}$ , $ s+m>\frac{3}{2}$. The system (\ref{-1.1}),(\ref{-1.2}),(\ref{-1.3}) with data $$(\psi_0,\phi_0,\phi_1) \in H^{s,2}({\bf R}^3) \times H^{m,2}({\bf R}^3) \times H^{m-1,2}({\bf R}^3) $$ has a unique global solution. More precisely, for any $T>0$ there exists a unique solution
$$(\psi,\phi,\phi_t) \in X^{s,\frac{1}{2}+}[0,T] \times \widetilde{X}^{m,\frac{1}{2}+}[0,T] \times \widetilde{X}^{m-1,\frac{1}{2}+}[0,T]. $$ 
This solution satisfies
$$ (\psi,\phi,\phi_t) \in C^0([0,T],H^{s,2}({\bf R}^3)\times H^{m,2}({\bf R}^3) \times H^{m-1,2}({\bf R}^3)) $$
Here $ \qquad \widetilde{X}^{m,\frac{1}{2}+}[0,T] := X_+^{m,\frac{1}{2}+}[0,T] + X_-^{m,\frac{1}{2}+}[0,T] $.
\end{theorem}
{\bf Remark:} The difference between the solution $(\psi,\phi,\phi_t)$ of the nonlinear problem (\ref{-1.1}),(\ref{-1.2}),(\ref{-1.3}) and the corresponding linear problem belongs to the space $H^{1,2}({\bf R}^3) \times H^{1,2}({\bf R}^3) \times L^2({\bf R}^3) $. This smoothing property (which is charateristic for the method) follows from Proposition \ref{Proposition 4.2} and \ref {Proposition 4.3}.

\end{document}